\documentclass[preprint,12pt]{elsarticle}

\usepackage{graphicx}
\usepackage{amsmath,amssymb,amsfonts, bm}
\usepackage{algorithmic}
\usepackage{textcomp}
\usepackage{xcolor}
\usepackage{hyperref}

\begin{document}

\begin{frontmatter}



\title{Social optimum of finite mean field games: existence and uniqueness of equilibrium solutions in  the finite horizon and stationary solutions in the infinite horizon\tnoteref{label6}}


\author[label1,label3]{Zijia Niu} 
\author[label1,label3]{Sanjin Huang}
\author[label1,label3]{Lu Ren}
\author[label2,label3,label4,label5]{Wang Yao \corref{cor1}}
\author[label1,label3,label4]{Xiao Zhang\corref{cor1}}

\affiliation[label1]{organization={The School of Mathematical Sciences},
            addressline={Beihang University}, 
            city={Beijing},
            postcode={100191}, 
            country={China}}
\affiliation[label2]{organization={The School of Artificial Intelligence},
            addressline={Beihang University}, 
            city={Beijing},
            postcode={100191}, 
            country={China}}
\affiliation[label3]{organization={Key Laboratory of Mathematics, Informatics and Behavioral Semantics},
            addressline={Ministry of Education}, 
            city={Beijing},
            postcode={100191}, 
            country={China}}
\affiliation[label4]{organization={Zhongguancun Laboratory},
            city={Beijing},
            postcode={100191}, 
            country={China}}
\affiliation[label5]{organization={Pengcheng Laboratory},
            city={Shenzhen},
            postcode={518000}, 
            state={Guangdong},
            country={China}}
\cortext[cor1]{Corresponding author}
 \tnotetext[label6]{This work was supported by the National Science and Technology Major Project under Grant 2022ZD0116401.}
\begin{abstract}
In this paper, we consider the social optimal problem of discrete time finite state space mean field games (referred to as finite mean field games \cite{HADIKHANLOO2019369}). Unlike the individual optimization of their own cost function in competitive models, in the problem we consider, individuals aim to optimize the social cost by finding a fixed point of the state distribution to achieve equilibrium in the mean field game. We provide a sufficient condition for the existence and uniqueness of the individual optimal strategies used to minimize the social cost. According to the definition of social optimum and the derived properties of social optimal cost,  the existence and uniqueness conditions of equilibrium solutions under initial-terminal value constraints in the finite horizon and the existence and uniqueness conditions of stable solutions in the infinite horizon are given. Finally, two examples that satisfy the conditions for the above solutions are provided.
\end{abstract}




\end{frontmatter}


\section{Introduction}
\label{section1}
Mean field games (MFGs) are developed by Lasry and Lions in a series of papers \cite{LASRY2006619, LASRY2006679, lasry2007mean}, and the related theory is proposed by Huang, Malhame, and Caines independently at almost the same time \cite{huang2006large, 4303232, 4434627, huang2007invariance}. MFGs are devoted to the analysis of differential games with a massive number of agents. These agents are homogeneous and indistinguishable. When the number of agents is large, the change of strategy of an individual has little impact on the overall system. It can be foreseen that it is unrealistic for an agent to make decisions by collecting detailed information from a massive number of other individuals. According to the MFGs theory, an agent only needs to consider the influence of the distribution of all individuals, thus simplifying the problem. MFGs have received widespread attention since their proposal, and related content can be found in references \cite{Gueant2011, Gomes2014, Achdou2020}.

Most of the literature related to MFGs considers continuous time games, and players are distributed on a continuum of states. Motivated both by the practical application requirements and numerical analysis questions that appear in the discretization of the continuous MFG system, discrete time finite state space mean field games have been proposed \cite{GOMES2010308}. This theory considers the transition of players between finite states, where players make decisions at discrete times and are influenced by the distribution of players in each state. Their strategy is the transition probability to each state. The Nash equilibrium of the game is obtained by solving the fixed point of a coupled system consisting of an individual cost optimization equation and a state distribution evolution equation. For simplicity, we adopt the nomenclature from the literature \cite{HADIKHANLOO2019369} to refer to discrete time finite state space mean field games as finite MFGs. Since its inception, finite MFGs have been applied in many fields, such as route choice \cite{8619448}, economic models \cite{BERNASCONI2023101974}, etc. Finite MFGs are competitive models in which each player optimizes their own cost function. In reality, there are often some scenarios where the social optimum needs to be considered. Social optimum refers to the cooperation of players to optimize the common social cost, which is the sum of individual costs \cite{9097879, Ma2023}.  The application scenarios involved include utility maximization in communication networks \cite{7869285} and the study of social welfare in economics \cite{Taylor2004}, etc.

In some literature on mean field Markov decision processes, there are discussions about the social optimum of systems in discrete time finite state space. These works are based on the mean field control, and the optimized objective function is the sum of the costs of individuals either in a finite horizon or with discount factors in the infinite time horizon. In a finite horizon, \cite{Gast2011} study the convergence of a Markove process of independent individuals evolving in a common environment to a simpler mean field model as the number of individuals tends to infinity, and also investigated the convergence rate. \cite{10.1214/21-AAP1713, Bauerle2023} consider the mean field Markov decision process with discounted reward in the infinite horizon. Both papers focus on the convergence problem of a Markov model with finite individuals to a mean field model when the number of individuals tends to infinity. In \cite{10.1214/21-AAP1713}, the authors deal with open-loop control and study the convergence rate. In \cite{Bauerle2023}, the authors prove the existence of optimal strategies, and also consider the problem of average reward. Similar related literature also includes \cite{6144708, Higuera-Chan2016}.

There is currently no research on social optimum in the literature related to finite MFGs. Through theoretical analysis of finite MFGs, it can be found that if the transition cost of an individual in state $i$ to state $j$ is only related to the decision of individuals in the state $i$, the optimal strategy corresponding to social optimum is the same as that corresponding to the competitive model. However, if the above transition cost is also related to the decisions of individuals in other states, the strategy selection under social optimum may differ from that of the competitive model. In this paper, we attempt to discuss the relevant issues of social optimum in finite MFGs, including the social optimum problem under initial-terminal value constraints in a finite horizon, as well as the social optimum problem in the infinite horizon. Our work is different from the above literature on mean field Markov processes. Firstly, in terms of the problem addressed, the above literature focuses on mean field control problems, while we focus on mean field game problems, where equilibrium solutions are obtained when the state distribution of individuals reaches a fixed point. Secondly, in terms of research content, the above literature mainly focuses on the convergence of models with finite individuals to mean field models when the number of individuals tends to infinity, while our work mainly focuses on the relevant properties of mean field game model, such as the existence and uniqueness conditions for the solution corresponding to the social optimum of finite MFGs. This is similar to the original paper on finite MFGs \cite{GOMES2010308}. In addition, \cite{10.1214/21-AAP1713, Bauerle2023} consider the problem with discounted rewards in the infinite horizon, while we study the stationary solutions in the infinite horizon without a discount factor. The main contributions of this paper are as follows:

(1) A sufficient condition for the existence and uniqueness of the optimal strategy corresponding to the social optimum of finite MFGs is given.

(2) The conditions for the existence and uniqueness of equilibrium solutions of the social optimum satisfying the initial-terminal constraints in a finite horizon are given.

(3) The conditions for the existence and uniqueness of stationary solutions of the social optimum in the infinite horizon are given.

The rest of the paper is organized as follows.
The social optima problem of finite MFGs is introduced in Section 2. In Section 3, we discuss the existence, uniqueness, and continuity of the optimal strategy for minimizing social cost. In Section 4, we provided the conditions for the existence and uniqueness of the equilibrium solution satisfying the initial-terminal constraints in a finite horizon. In Section 5, we prove some propositions, lemmas, and theorems related to the existence and uniqueness of stationary solutions in the infinite horizon. In Section 6, two examples that correspond to solutions of finite horizon and infinite horizon are provided. Finally, the paper is concluded in Section 7.

\section{The Social optima problem of finite MFGs}
\label{section2}
Consider a system with a large number of homogeneous agents that make decisions at discrete times to transition between several finite states. The set of finite states is denoted as $\mathcal{S}=\{1, 2, ... ,s\}$ and the set of decision times is denoted as $\mathcal{N}=\{0, 1, ... ,N-1\}$. Define the probability simplex $\mathbb{P}=\left\{\left[p_1, ...,p_d\right]|p_j\geq0\;\forall j, \sum_jp_j=1\right\}$, then the distribution of agents among states at time $n$ can be supposed as $\boldsymbol{m}^n=\left[m^n_1,...,m^n_s\right]\in\mathbb{P}$, and the strategy of agents at time $n$, also known as the state transition probability matrix, can be represented as $\boldsymbol{P}^n=[P^n_{ij}]_{i,j\in\mathcal{S}}\in\mathbb{P}^s$. We can recursively write the evolution equations of the agents' distribution at adjacent times. For any $n\in\mathcal{N}$,
\begin{equation}\label{deq}
m^{n+1}_j=\sum_{i\in\mathcal{S}}m^n_iP^n_{ij}.
\end{equation}
The transition cost of an agent from state $i$ to state $j$ at time $n$ is defined as $c_{ij}(\boldsymbol{m}^n, \boldsymbol{P}^n)$. Set the terminal cost to $\boldsymbol{G}^N\in\mathbb{R}^s$, then the social cost of the system within a given time horizon $\mathcal{N}\cup\{N\}$ is represented as
\begin{equation}\label{socialcost}
H^0(\boldsymbol{m}, \boldsymbol{P}^{\{n\}})=\sum_{n=0}^{N-1}\sum_{i,j\in\mathcal{S}}c_{ij}(\boldsymbol{m}^n, \boldsymbol{P}^n)m^n_iP^n_{ij}+\sum_jG^N_jm^N_j.
\end{equation}

The similar recursive form of social cost at any time $n\in\mathcal{N\cup\{N\}}$ can be derived as follows. Let $\hat{H}^n=\mathop{min}_{\boldsymbol{P}\in\mathbb{P}^s}H^n=\sum_{i\in\mathcal{S}}U^n_im^n_i$ be the optimal social cost at time $n\in\mathcal{N}\cup\{N\}$, where $U^n_i$ is the optimal cost of an agent in state $i$ at time $n$. At time $N$, $U^N_i=G^N_i$, at time $n\in\mathcal{N}$, $U^n_i$ can be defined as
\begin{equation}\label{indicost}
U^n_i=\sum_{j\in\mathcal{S}}\left(c_{ij}(\boldsymbol{m}^n, \hat{\boldsymbol{P}}^n)\hat{P}^n_{ij}+U^{n+1}_j\hat{P}^n_{ij}\right),
\end{equation}
where $\hat{\boldsymbol{P}}^n=\mathop{argmin}_{\boldsymbol{P}\in\mathbb{P}^s}H^n$. It can be seen that the ``optimal" of $U^n_i$ here does not mean that $U^n_i$ is taken to the minimum, but rather the individual cost of state $i$, which is pushed back by the strategy corresponding to minimizing the social cost. Then, based on the Dynamic Programming Principle, for any $n\in\mathcal{N}$,
\begin{equation}\label{recusc}
\begin{split}
\mathop{min}_{\boldsymbol{P}\in\mathbb{P}^s}H^n(\boldsymbol{m}^n, \boldsymbol{P}^n, \boldsymbol{U}^{n+1})&=\mathop{min}_{\boldsymbol{P}\in\mathbb{P}^s}\sum_{i,j\in\mathcal{S}}c_{ij}(\boldsymbol{m}^n, \boldsymbol{P}^n)m^n_iP^n_{ij}+\hat{H}^{n+1} \\
                                                                                                                                &=\mathop{min}_{\boldsymbol{P}\in\mathbb{P}^s}\sum_{i,j\in\mathcal{S}}c_{ij}(\boldsymbol{m}^n, \boldsymbol{P}^n)m^n_iP^n_{ij}+\sum_{j\in\mathcal{S}}U^{n+1}_jm^{n+1}_j \\
                                                                                                                                &=\mathop{min}_{\boldsymbol{P}\in\mathbb{P}^s}\sum_{i,j\in\mathcal{S}}c_{ij}(\boldsymbol{m}^n, \boldsymbol{P}^n)m^n_iP^n_{ij}+\sum_{i, j\in\mathcal{S}}U^{n+1}_jm^n_iP^n_{ij}.
\end{split}
\end{equation}

The definition of the optimal strategy $\hat{\boldsymbol{P}}^n$ for minimizing social cost is given below.

\textbf{Definition 2.1.} For any $n\in\mathcal{N}$, fix a state distribution $\boldsymbol{m}^n\in\mathbb{P}$ and a cost vector $\boldsymbol{U}^{n+1}\in\mathbb{R}^s$, if for any $i\in\mathcal{S}$, $\boldsymbol{q}\in\mathbb{P}$,
\begin{equation}\label{os}
H^n(\boldsymbol{m}^n, \hat{\boldsymbol{P}}^n, \boldsymbol{U}^{n+1})\leq H^n(\boldsymbol{m}^n, \mathcal{P}(\hat{\boldsymbol{P}}^n, \boldsymbol{q}, i), \boldsymbol{U}^{n+1}),
\end{equation}
where $\mathcal{P}(\hat{\boldsymbol{P}}^n, \boldsymbol{q}, i)\in\mathbb{P}^s$ represents the strategy matrix formed by replacing the $i$-th row of $\hat{\boldsymbol{P}}^n$ with $\boldsymbol{q}$. Then $\hat{\boldsymbol{P}}^n$ is an optimal strategy corresponding to the minimum social cost of finite MFGs at time $n$. Correspondingly, in the decision time set $\mathcal{N}$, the strategy sequences of all agents obtained based on $(\hat{\boldsymbol{P}}^n)_{n\in\mathcal{N}}$ constitute the Nash equilibrium of the game.

The social optima problem of finite MFGs can be defined as follows:

\textbf{Definiton 2.2.} (P1)  Given the initial state distribution $\boldsymbol{m}^0=\tilde{\boldsymbol{m}}$ and the terminal cost $\boldsymbol{U}^N=\boldsymbol{G}^N$, find the optimal strategy sequence $\left(\hat{\boldsymbol{P}}^n\right)_{n\in\mathcal{N}}$ that satisfies the following constraints
\begin{equation}\label{oscs}
\left\{
\begin{split}
&\hat{\boldsymbol{P}}^n=\mathop{argmin}_{\boldsymbol{P}^n\in\mathbb{P}^s} H^n(\boldsymbol{m}^n, \boldsymbol{P}^n, \boldsymbol{U}^{n+1}) \\
&\sum_{j\in\mathcal{S}}\hat{P}^n_{ij}=1, \quad i\in\mathcal{S} \\
&\hat{P}^n_{ij}\geq0, \quad i,j\in\mathcal{S},
\end{split}
\right.
\end{equation}
such that a sequence of pairs of s-dimensional vectors $\left(\boldsymbol{m}^n, \boldsymbol{U}^n\right)_{n\in\mathcal{N}\cup\{N\}}$ satisfies 
\begin{equation}\label{fmfgso}
\left\{
\begin{split}
&U^n_i=\sum_j(c_{ij}(\boldsymbol{m}^n, \hat{\boldsymbol{P}}^n)\hat{P}^n_{ij}+U^{n+1}_j\hat{P}^n_{ij}) \\
&m^{n+1}_j=\sum_im^n_i\hat{P}^n_{ij}.
\end{split}
\right.
\end{equation}

\section{Properties of the optimal strategy}
\label{section3}
In this section, we discuss the relevant properties of the optimal strategy $\hat{\boldsymbol{P}}^n$ that minimizes the social cost $H^n(\boldsymbol{m}^n, \boldsymbol{P}^n, \boldsymbol{U}^{n+1})$ at any time $n\in\mathcal{N}$. The assumptions involved are listed below:

\textbf{(A1)} For any $\boldsymbol{m}\in\mathbb{P}, \boldsymbol{P}\in\mathbb{P}^s, i\in\mathcal{S}, j\in\mathcal{S}$, $c_{ij}(\boldsymbol{m}, \boldsymbol{P})$ is a differentiable convex function with respect to $\boldsymbol{P}$. And for any $\boldsymbol{P}^1, \boldsymbol{P}^2\in\mathbb{P}^s$, $\boldsymbol{P}^1\neq\boldsymbol{P}^2$, the following inequality holds
\begin{equation} \label{as1}
\sum_{i, j\in\mathcal{S}}\left(c_{ij}(\boldsymbol{m}, \boldsymbol{P}^1)-c_{ij}(\boldsymbol{m}, \boldsymbol{P}^2)\right)\left(P^1_{ij}-P^2_{ij}\right)m_i\geq0.
\end{equation}

\textbf{(A2)}  $H(\boldsymbol{m}, \boldsymbol{P}, \boldsymbol{U}) : \mathbb{P}\times\mathbb{P}^s\times\mathbb{R}^s\rightarrow\mathbb{R}$ is a continuous function.

\subsection{The existence and uniqueness of the optimal strategy}
\label{section3.1}
Consider the following optimization problem (P2):
\begin{equation}\nonumber
\begin{split}
\mathop{min}_{\boldsymbol{P}^n\in\mathbb{P}^s}\quad&H^n(\boldsymbol{m}^n, \boldsymbol{P}^n, \boldsymbol{U}^{n+1}) \\
s.t. \quad&\sum_{j\in\mathcal{S}}P^n_{ij}=1, \quad i\in\mathcal{S} \\
      \quad&P^n_{ij}\geq0, \quad i,j\in\mathcal{S}. 
\end{split}
\end{equation}
The following theorem provides a sufficient condition for the existence and uniqueness of the optimal strategy for the problem (P2).

\textbf{Theorem 3.1.} For any $\boldsymbol{m}^n\in\mathbb{P}, \boldsymbol{U}^{n+1}\in\mathbb{R}^s, n\in\mathcal{N}$ that hold assumption (A1), there exists an optimal strategy is solution to problem (P2), if only the greater-than sign in \eqref{as1} holds, the optimal strategy is unique.

\textbf{Proof.} \textbf{Existence} According to (A1), $c_{ij}(\boldsymbol{m}^n, \boldsymbol{P}^n)$ is a differentiable convex function with respect to $\boldsymbol{P}^n$, then for any $\boldsymbol{P}^{n,1}, \boldsymbol{P}^{n,2}\in\mathbb{P}^s$, any $\alpha\in(0, 1)$, the following inequality holds
\begin{equation}\label{oseu1}
c_{ij}\left(\boldsymbol{m}^n, \alpha\boldsymbol{P}^{n, 1}+(1-\alpha)\boldsymbol{P}^{n, 2}\right)\leq\alpha c_{ij}(\boldsymbol{m}^n, \boldsymbol{P}^{n, 1})+(1-\alpha)c_{ij}(\boldsymbol{m}^n, \boldsymbol{P}^{n, 2}).
\end{equation}
Let $g(\boldsymbol{P}^{n})=\sum_{i, j\in\mathcal{S}}c_{ij}(\boldsymbol{m}^n, \boldsymbol{P}^{n})m^n_iP^n_{ij}$, then
\begin{equation}\label{oseu2}
\begin{split}
g(\alpha\boldsymbol{P}^{n,1}+(1-\alpha)\boldsymbol{P}^{n,2})=&\sum_{i, j\in\mathcal{S}}c_{ij}\left(\boldsymbol{m}^n, \alpha\boldsymbol{P}^{n, 1}\!+\!(1\!-\!\alpha)\boldsymbol{P}^{n,2}\right)m^n_i\left(\alpha P^{n, 1}_{ij}\!+\!(1\!-\!\alpha)P^{n, 2}_{ij}\right) \\                                                                                                             
                                                                                                              \leq&\sum_{i, j\in\mathcal{S}}m^n_i\left[\alpha^2P^{n, 1}_{ij}c_{ij}(\boldsymbol{m}^n, \boldsymbol{P}^{n, 1})\!+\!\alpha(1\!-\!\alpha)P^{n, 1}_{ij}c_{ij}(\boldsymbol{m}^n, \boldsymbol{P}^{n,2})+\right. \\
                                                                                                                    &\left. \alpha(1\!-\!\alpha)P^{n, 2}_{ij}c_{ij}(\boldsymbol{m}^n, \boldsymbol{P}^{n, 1})\!+\!(1\!-\!\alpha)^2P^{n, 2}_{ij}c_{ij}(\boldsymbol{m}^n, \boldsymbol{P}^{n, 2})\right].
\end{split}
\end{equation}
Let
\begin{equation}\nonumber
\begin{split}
h(\boldsymbol{P}^{n, 1}, \boldsymbol{P}^{n, 2})=&\sum_{i, j\in\mathcal{S}}m^n_i\left[\alpha^2P^{n, 1}_{ij}c_{ij}(\boldsymbol{m}^n, \boldsymbol{P}^{n, 1})+\alpha(1-\alpha)P^{n, 1}_{ij}c_{ij}(\boldsymbol{m}^n, \boldsymbol{P}^{n, 2})+\right. \\
                                                                                        &\left. \alpha(1-\alpha)P^{n, 2}_{ij}c_{ij}(\boldsymbol{m}^n, \boldsymbol{P}^{n, 1})+(1-\alpha)^2P^{n, 2}_{ij}c_{ij}(\boldsymbol{m}^n, \boldsymbol{P}^{n, 2})\right],
\end{split}
\end{equation}
we have
\begin{equation}\label{oseu3}
\begin{split}
  &\alpha g(\boldsymbol{P}^{n, 1})+(1-\alpha)g(\boldsymbol{P}^{n, 2})-h(\boldsymbol{P}^{n, 1}, \boldsymbol{P}^{n, 2}) \\
=&\sum_{i, j\in\mathcal{S}}m^n_i\left[\alpha(1-\alpha)\left(c_{ij}(\boldsymbol{m}^n, \boldsymbol{P}^{n, 1})(P^{n, 1}_{ij}-P^{n, 2}_{ij})+c_{ij}(\boldsymbol{m}^n, \boldsymbol{P}^{n, 2})(P^{n, 2}_{ij}-P^{n, 1}_{ij})\right)\right]\\
=&\alpha(1-\alpha)\sum_{i, j\in\mathcal{S}}\left(c_{ij}(\boldsymbol{m}^n, \boldsymbol{P}^{n,1})-c_{ij}(\boldsymbol{m}^n, \boldsymbol{P}^{n,2})\right)\left(P^{n, 1}_{ij}-P^{n, 2}_{ij}\right)m^n_i \\
\geq&0.
\end{split}
\end{equation}
According to assumption (A1), the greater-than or equal sign in the above equation holds.

For any $\boldsymbol{P}^{n,1}, \boldsymbol{P}^{n,2}\in\mathbb{P}^s$, any $\alpha\in(0, 1)$, the following inequality can be obtained from equations \eqref{oseu2} and \eqref{oseu3}
\begin{equation}\label{oseu4}
g\left(\alpha\boldsymbol{P}^{n,1}+(1-\alpha)\boldsymbol{P}^{n,2}\right)\leq\alpha g(\boldsymbol{P}^{n, 1})+(1-\alpha)g(\boldsymbol{P}^{n, 2}).
\end{equation}
Then $g(\boldsymbol{P}^n)$ is a convex function with respect to $\boldsymbol{P}^n$. According to \eqref{recusc}, $H^n(\boldsymbol{m}^n,\boldsymbol{P}^n, \boldsymbol{U}^{n+1}), n\in\mathcal{N}$ is a convex function with respect to $\boldsymbol{P}^n$. We focus on the feasible set $\mathbb{P}^s$ of $\boldsymbol{P}^n$. Since the set $\mathbb{P}^s$ is a union of $s$ probability simplexes, it is a non-empty compact set. Since $H^n(\boldsymbol{m}^n,\boldsymbol{P}^n, \boldsymbol{U}^{n+1})$ is a convex function with respect to $\boldsymbol{P}^n$, and the constraints in problem (P2) are linear constraints about $\boldsymbol{P}^n$, problem (P2) is a convex optimization. Because the feasible set $\mathbb{P}^s$ is a non-empty compact set, there exists $\hat{\boldsymbol{P}}^n$ as the optimal strategy for the convex optimization problem (P2).

\textbf{Uniqueness} \textbf{(Proof 1)} If only the greater-than sign in \eqref{as1} holds, then according to \eqref{oseu2} and \eqref{oseu3}, \eqref{oseu4} becomes
\begin{equation}\label{oseu5}
g\left(\alpha\boldsymbol{P}^{n,1}+(1-\alpha)\boldsymbol{P}^{n,2}\right)<\alpha g(\boldsymbol{P}^{n, 1})+(1-\alpha)g(\boldsymbol{P}^{n, 2}).
\end{equation}
Then $g(\boldsymbol{P}^n)$ is a strictly convex function with respect to $\boldsymbol{P}^n$, and $H^n$ is a strictly convex function with respect to $\boldsymbol{P}^n$. Therefore, there exists a unique optimal strategy for the convex optimization problem (P2) on the feasible set $\mathbb{P}^s$.

In addition to the above methods, we can also prove that the optimal strategy is unique when only the greater-than sign in \eqref{as1} holds, starting from the composition of the social cost.

\textbf{Uniqueness} \textbf{(Proof 2)} If both $\hat{\boldsymbol{P}}^{n,1}$ and $\hat{\boldsymbol{P}}^{n,2}$ are optimal strategies, then for each component $\hat{P}^n_{ij}$, they need to satisfy the follow Karush-Kuhn-Tucker (KKT) conditions:
\begin{subequations}
\begin{align}
&\frac{\partial H^n(\boldsymbol{m}^n, \hat{\boldsymbol{P}}^{n, k}, \boldsymbol{U}^{n+1})}{\partial\hat{P}^n_{ij}}-\lambda^{n, k}_{i}-\mu^{n, k}_{ij}=0, \label{oseu6a} \\
&\sum_{j\in\mathcal{S}}\hat{P}^{n, k}_{ij}-1=0, \label{oseu6b}\\
&\hat{P}^{n, k}_{ij}\geq0, \label{oseu6c} \\
&\mu^{n, k}_{ij}\geq0, \label{oseu6d}\\
&\mu^{n, k}_{ij}\hat{P}^{n, k}_{ij}=0, \label{oseu6e}
\end{align}
\end{subequations}
where $k\in\{1, 2\}$. $\lambda^n_i$ is the Lagrange multiplier associated with $\sum_{j\in\mathcal{S}}P^n_{ij}=1$, and $\mu^n_{ij}\geq0$ corresponds to the inequality constraint $P^n_{ij}\geq0$. According to \eqref{recusc} we have
\begin{equation}\label{oseu7}
\frac{\partial H^n(\boldsymbol{m}^n, \hat{\boldsymbol{P}}^n, \boldsymbol{U}^{n+1})}{\partial\hat{P}^n_{ij}}=\sum_{p, q\in\mathcal{S}}\frac{\partial c_{pq}(\boldsymbol{m}^n, \hat{\boldsymbol{P}}^n)}{\partial\hat{P}^n_{ij}}m^n_p\hat{P}^n_{pq}+c_{ij}(\boldsymbol{m}^n, \hat{\boldsymbol{P}^n})m^n_i+U^{n+1}_jm^n_i.
\end{equation}
Due to the arbitrary selection of component $\hat{P}^n_{ij}$, for all $i, j\in\mathcal{S}$, from \eqref{oseu6a} and \eqref{oseu7} we can obtain that 
\begin{equation}\label{oseu8}
\begin{split}
&\sum_{i, j\in\mathcal{S}}\left(\hat{P}^{n, 1}_{ij}\!-\!\hat{P}^{n, 2}_{ij}\right)\left[\sum_{p, q\in\mathcal{S}}\left(\frac{\partial c_{pq}(\boldsymbol{m}^n, \hat{\boldsymbol{P}}^{n, 1})}{\partial\hat{P}^{n, 1}_{ij}}m^n_p\hat{P}^{n, 1}_{pq}\!-\!\frac{\partial c_{pq}(\boldsymbol{m}^n, \hat{\boldsymbol{P}}^{n, 2})}{\partial\hat{P}^{n, 2}_{ij}}m^n_p\hat{P}^{n, 2}_{pq}\right)\!+\right. \\
&\left. c_{ij}(\boldsymbol{m}^n, \hat{\boldsymbol{P}}^{n, 1})m^n_i-c_{ij}(\boldsymbol{m}^n, \hat{\boldsymbol{P}}^{n, 2})m^n_i-\lambda^{n, 1}_{i}-\mu^{n, 1}_{ij}+\lambda^{n, 2}_{i}+\mu^{n, 2}_{ij}\right]=0.
\end{split}
\end{equation}
Consider
\begin{equation}\label{oseu9}
\begin{split}
&\sum_{i, j\in\mathcal{S}}\left(\hat{P}^{n, 1}_{ij}\!-\!\hat{P}^{n, 2}_{ij}\right)\left[\sum_{p, q\in\mathcal{S}}\left(\frac{\partial c_{pq}(\boldsymbol{m}^n, \hat{\boldsymbol{P}}^{n, 1})}{\partial\hat{P}^{n, 1}_{ij}}m^n_p\hat{P}^{n, 1}_{pq}\!-\!\frac{\partial c_{pq}(\boldsymbol{m}^n, \hat{\boldsymbol{P}}^{n, 2})}{\partial\hat{P}^{n, 2}_{ij}}m^n_p\hat{P}^{n, 2}_{pq}\right)\!+\right. \\
&\left.c_{ij}(\boldsymbol{m}^n, \hat{\boldsymbol{P}}^{n, 1})m^n_i-c_{ij}(\boldsymbol{m}^n, \hat{\boldsymbol{P}}^{n, 2})m^n_i\right] \\
=&\sum_{p,q\in\mathcal{S}}\left[-\sum_{i,j\in\mathcal{S}}\frac{\partial c_{pq}(\boldsymbol{m}^n, \hat{\boldsymbol{P}}^{n, 1})}{\partial\hat{P}^{n, 1}_{ij}}\left(\hat{P}^{n, 2}_{ij}-\hat{P}^{n, 1}_{ij}\right)m^n_p\hat{P}^{n, 1}_{pq}-\right. \\
&\quad\quad\left. \sum_{i,j\in\mathcal{S}}\frac{\partial c_{pq}(\boldsymbol{m}^n, \hat{\boldsymbol{P}}^{n, 2})}{\partial\hat{P}^{n, 2}_{ij}}\left(\hat{P}^{n, 1}_{ij}-\hat{P}^{n, 2}_{ij}\right)m^n_p\hat{P}^{n, 2}_{pq}\right] \\
&+\sum_{p,q\in\mathcal{S}}\left(c_{pq}(\boldsymbol{m}^n, \hat{\boldsymbol{P}}^{n, 1})-c_{pq}(\boldsymbol{m}^n, \hat{\boldsymbol{P}}^{n, 2})\right)\left(\hat{P}^{n, 1}_{pq}-\hat{P}^{n, 2}_{pq}\right)m^n_p,
\end{split}
\end{equation}
According to assumption (A1), $c_{ij}(\boldsymbol{m}^n, \boldsymbol{P}^n)$ is a convex function with respect to $\boldsymbol{P}^n$, then for any $\boldsymbol{P}^n, \tilde{\boldsymbol{P}}^n\in\mathbb{P}^s$, the following inequality holds
\begin{equation}\label{oseu10}
c_{pq}(\boldsymbol{m}^n, \boldsymbol{P}^n)-c_{pq}(\boldsymbol{m}^n, \tilde{\boldsymbol{P}}^n)-\sum_{i,j\in\mathcal{S}}\frac{\partial c_{pq}(\boldsymbol{m}^n, \tilde{\boldsymbol{P}}^n)}{\partial\tilde{P}^n_{ij}}(P^n_{ij}-\tilde{P}^n_{ij})\geq0,
\end{equation}
then \eqref{oseu9} satisfies
\begin{equation}\label{oseu11}
\begin{split}
\eqref{oseu9}\geq&\sum_{p, q\in\mathcal{S}}\left[\left(c_{pq}(\boldsymbol{m}^n,\hat{\boldsymbol{P}}^{n, 1})-c_{pq}(\boldsymbol{m}^n,\hat{\boldsymbol{P}}^{n, 2})\right)m^n_p\hat{P}^{n,1}_{pq}-\right. \\
&\quad\quad\left. \left(c_{pq}(\boldsymbol{m}^n,\hat{\boldsymbol{P}}^{n, 1})\!-\!c_{pq}(\boldsymbol{m}^n,\hat{\boldsymbol{P}}^{n, 2})\right)m^n_p\hat{P}^{n,2}_{pq}\right] \\
&+\sum_{p, q\in\mathcal{S}}\left[\left(c_{pq}(\boldsymbol{m}^n,\hat{\boldsymbol{P}}^{n, 1})-c_{pq}(\boldsymbol{m}^n,\hat{\boldsymbol{P}}^{n, 2})\right)\left(\hat{P}^{n,1}_{pq}-\hat{P}^{n,2}_{pq}\right)m^n_p\right] \\
=&2\sum_{p, q\in\mathcal{S}}\left[\left(c_{pq}(\boldsymbol{m}^n,\hat{\boldsymbol{P}}^{n, 1})-c_{pq}(\boldsymbol{m}^n,\hat{\boldsymbol{P}}^{n, 2})\right)\left(\hat{P}^{n,1}_{pq}-\hat{P}^{n,2}_{pq}\right)m^n_p\right] \\
>&0.
\end{split}
\end{equation}
From \eqref{oseu8} and \eqref{oseu11} we can conclude that
\begin{equation}\label{oseu12}
\sum_{i, j\in\mathcal{S}}\left(\hat{P}^{n, 1}_{ij}-\hat{P}^{n, 2}_{ij}\right)\left(-\lambda^{n, 1}_{i}-\mu^{n, 1}_{ij}+\lambda^{n, 2}_{i}+\mu^{n, 2}_{ij}\right)<0.
\end{equation}
Using  $\mu^{n, k}_{ij}\hat{P}^{n, k}_{ij}=0$, there is
\begin{equation}\label{oseu13}
\sum_{i, j\in\mathcal{S}}\left(-\hat{P}^{n, 1}_{ij}\lambda^{n, 1}_i+\hat{P}^{n, 1}_{ij}\lambda^{n, 2}_i+\hat{P}^{n, 1}_{ij}\mu^{n, 2}_{ij}+\hat{P}^{n, 2}_{ij}\lambda^{n, 1}_i+\hat{P}^{n, 2}_{ij}\mu^{n, 1}_{ij}-\hat{P}^{n, 2}_{ij}\lambda^{n, 2}_i \right)<0.
\end{equation}
Since $\sum_{j\in\mathcal{S}}\hat{P}^{n, k}_{ij}\lambda^{n, l}_i=\lambda^{n, l}_i$, where $k,l\in\{1,2\}$, \eqref{oseu13} can be transformed into
\begin{equation}\label{oseu14}
\sum_{i, j\in\mathcal{S}}\left(\hat{P}^{n, 1}_{ij}\mu^{n, 2}_{ij}+\hat{P}^{n, 2}_{ij}\mu^{n, 1}_{ij}\right)<0.
\end{equation}
According to $\hat{P}^{n, k}_{ij},\mu^{n, k}_{ij}\geq0$, we obtain a contradiction.
$\hfill\square$

Since Theorem 3.1 only provides a sufficient condition for the existence and uniqueness of the optimal strategy corresponding to the social optimum, the optimal strategy can also be obtained through other assumptions. Therefore, the following assumption is provided for ease of use in subsequent content.

\textbf{(A3)} For any $n\in\mathcal{N}$, the optimization problem (P2) has a unique optimal strategy $\hat{\boldsymbol{P}}^n$.

\subsection{The continuity of the optimal strategy}
The proof that the optimal strategy $\hat{\boldsymbol{P}}^n$ is a continuous function of $\boldsymbol{m}^n$ and $\boldsymbol{U}^{n+1}$ is given as follows.

\textbf{Theorem 3.2.} Suppose assumptions (A2) and (A3) hold, then for any $n\in\mathcal{N}$, the optimal strategy $\hat{\boldsymbol{P}}^n$ is a continuous function of $\boldsymbol{m}^n$ and $\boldsymbol{U}^{n+1}$.

\textbf{Proof.} Consider sequences $\boldsymbol{m}^n_{(k)}\rightarrow\boldsymbol{m}^n$ and $\boldsymbol{U}^{n+1}_{(k)}\rightarrow\boldsymbol{U}^{n+1}$, for an optimal strategy $\tilde{\boldsymbol{P}}^n$, $H^n(\boldsymbol{m}^n, \tilde{\boldsymbol{P}}^n, \boldsymbol{U}^{n+1})$ is a continuous function, then there exists a sequence $\tilde{\boldsymbol{P}}^n_{(k)}=\tilde{\boldsymbol{P}}^n(\boldsymbol{m}^n_{(k)}, \boldsymbol{U}^{n+1}_{(k)})$ that converges to the optimal strategy $\tilde{\boldsymbol{P}}^n$. According to  (A3), the optimal strategy is unique, that is, $\tilde{\boldsymbol{P}}^n=\hat{\boldsymbol{P}}^n$, then $\tilde{\boldsymbol{P}}^n_{(k)}\rightarrow\hat{\boldsymbol{P}}^n$, $\hat{\boldsymbol{P}}^n$ is continuous with respect to $\boldsymbol{m}^n$ and $\boldsymbol{U}^{n+1}$.
$\hfill\square$

\section{Solutions to the social optima problem of finite MFGs with initial-terminal value}
\label{section4}
In this section, we discuss the existence and uniqueness of solutions to the social optima problem of finite MFGs with given initial-terminal values (problem (P1)). Compared to the competitive model in the original paper on finite MFGs \cite{GOMES2010308}, the uniqueness condition for the solutions to the social optima problem is extended. In the competitive model, the transition cost $c_{i*}(\boldsymbol{m}^n, \boldsymbol{P}^n_i)$ at state $i$ only related to the strategy  $\boldsymbol{P}^n_i$ of players in state $i$. In the social optima model, $c_{i*}(\boldsymbol{m}^n, \boldsymbol{P}^n)$ can be influenced by $\boldsymbol{P}^n$ which includes the strategies of players in other states.

\subsection{Existence}
\label{section4.1}
\textbf{Theorem 4.1.} Given any initial distribution $\boldsymbol{m}^0=\tilde{\boldsymbol{m}}\in\mathbb{P}$ and terminal cost $\boldsymbol{U}^N=\boldsymbol{G}^N\in\mathbb{R}^s$ that make assumption (A2) and (A3) hold, there exists a sequence of pairs of s-dimensional vectors
\begin{equation}\nonumber
\left(\hat{\boldsymbol{m}}^n, \hat{\boldsymbol{U}}^n\right)_{n\in\mathcal{N}\cup\{N\}}
\end{equation}
that is a solution to the social optima problem of finite MFGs with given initial-terminal values (problem (P1)).

\textbf{Proof.} Assuming a state distribution sequence $\left(\boldsymbol{m}^{n,[1]}\right)_{n\in\mathcal{N}\cup\{N\}}\in\mathbb{P}^{N+1}$ is given, where the initial distribution is $\boldsymbol{m}^{0,[1]}=\tilde{\boldsymbol{m}}$. Then we can obtain the optimal social cost sequence $\left(H^{n}\right)_{n\in\mathcal{N}\cup\{N\}}$ and corresponding unique optimal strategy sequence $\left(\hat{\boldsymbol{P}}^{n}\right)_{n\in\mathcal{N}}$ based on the following two equations, as well as the individual cost sequence $\left(\boldsymbol{U}^{n}\right)_{n\in\mathcal{N}\cup\{N\}}$ of each state.
\begin{equation}\label{itfmfgsoexist1}
H^n(\boldsymbol{m}^{n,[1]}, \hat{\boldsymbol{P}}^n, \boldsymbol{U}^{n+1})=\sum_{i,j\in\mathcal{S}}c_{ij}(\boldsymbol{m}^{n,[1]}, \hat{\boldsymbol{P}}^n)m^{n,[1]}_i\hat{P}^n_{ij}+\sum_{i,j\in\mathcal{S}}U^{n+1}_jm^{n,[1]}_i\hat{P}^n_{ij},
\end{equation}
\begin{equation}\label{itfmfgsoexist2}
U^n_i=\sum_{j\in\mathcal{S}}\left(c_{ij}(\boldsymbol{m}^{n,[1]},\hat{\boldsymbol{P}}^n)\hat{P}^n_{ij}+U^{n+1}_j\hat{P}^n_{ij}\right).
\end{equation}
Let $\boldsymbol{m}^{0,[2]}=\tilde{\boldsymbol{m}}$ and use
\begin{equation}\label{itfmfgsoexist3}
\boldsymbol{m}^{n+1,[2]}_j=\sum_{i\in\mathcal{S}}m^{n,[2]}_i\hat{P}^n_{ij},
\end{equation}
we can obtain the second state distribution sequence $\{\boldsymbol{m}^{n,[2]}\}_{n\in\mathcal{N}\cup\{N\}}\in\mathbb{P}^{N+1}$. According to the continuity of \eqref{itfmfgsoexist1}-\eqref{itfmfgsoexist3} (Theorem 3.2), it can be seen that the above process defines a continuous mapping from $\mathbb{P}^{N+1}$ to $\mathbb{P}^{N+1}$, which maps $\left(\boldsymbol{m}^{n,[1]}\right)_{n\in\mathcal{N}\cup\{N\}}$ to $\left(\boldsymbol{m}^{n,[2]}\right)_{n\in\mathcal{N}\cup\{N\}}$. According to Brower's fixed point theorem, there exists a solution $\left(\hat{\boldsymbol{m}}^n, \hat{\boldsymbol{U}}^n\right)_{n\in\mathcal{N}\cup\{N\}}$ corresponding to the social optima problem of finite MFGs with given initial-terminal values (problem (P1)).
$\hfill\square$

\subsection{Uniqueness}
\label{section4.2}
At any time $n\in\mathcal{N}$, if the optimal strategy $\hat{\boldsymbol{P}}^n$ exists and is unique, we can define the optimal social cost as a function of individual cost in each state at time $n+1$, that is
\begin{equation}\label{itfmfgsouniq1}
\begin{split}
\Phi_{\boldsymbol{m}^n}(\boldsymbol{U}^{n+1})=&\mathop{min}_{\boldsymbol{P}^n}H^n(\boldsymbol{m}^n, \boldsymbol{P}^n, \boldsymbol{U}^{n+1})\\
=&\sum_{i,j\in\mathcal{S}}c_{ij}(\boldsymbol{m}^n, \hat{\boldsymbol{P}}^n)m^n_i\hat{P}^n_{ij}+\sum_{i,j\in\mathcal{S}}U^{n+1}_jm^n_i\hat{P}^n_{ij}.
\end{split}
\end{equation}
Similarly, based on the optimal strategy $\hat{\boldsymbol{P}}^n$, the recursive function of individual costs corresponding to the social optimum can be defined as
\begin{equation}\label{itfmfgsouniq2}
\Gamma_{\boldsymbol{m}^n}(\boldsymbol{U}^{n+1}\!)\!\!=\!\!\left[\sum_{j\in\mathcal{S}}\!\left(c_{1j}(\!\boldsymbol{m}^n, \!\hat{\boldsymbol{P}}^n\!)\hat{P}^n_{1j}\!+\!U^{n+1}_j\!\hat{P}^n_{1j}\right)\!,\!...,\!\!\sum_{j\in\mathcal{S}}\!\left(c_{sj}(\!\boldsymbol{m}^n, \!\hat{\boldsymbol{P}}^n\!)\hat{P}^n_{sj}\!+\!U^{n+1}_j\!\hat{P}^n_{sj}\right)\!\right]\!,
\end{equation}
which belongs to $\mathbb{R}^s$. The recursive function of the state distribution can be defined as
\begin{equation}\label{itfmfgsouniq3}
\Theta_{\boldsymbol{U}^{n+1}}(\boldsymbol{m}^n)=\boldsymbol{m}^n\hat{\boldsymbol{P}}^n.
\end{equation}
Before proving the uniqueness of the solution, a lemma is provided.

\textbf{Lemma 4.1.} For any $\boldsymbol{m}\in\mathbb{P}$ and $\boldsymbol{U}^1, \boldsymbol{U}^2\in\mathbb{R}^s$ that hold assumption (A3), the optimal social cost $\Phi_{\boldsymbol{m}}(\boldsymbol{U})$ satisfies the following inequality
\begin{equation}\label{itfmfgsouniq4}
\Phi_{\boldsymbol{m}}(\boldsymbol{U}^2)-\Phi_{\boldsymbol{m}}(\boldsymbol{U}^1)-\left(\boldsymbol{U}^2-\boldsymbol{U}^1\right)\cdot\Theta_{\boldsymbol{U}^1}(\boldsymbol{m})\leq0.
\end{equation}

\textbf{Proof.} Assuming $\hat{\boldsymbol{P}}^1, \hat{\boldsymbol{P}}^2\in\mathbb{P}^s$ are the optimal strategies correspongding to $\Phi_{\boldsymbol{m}}(\boldsymbol{U}^1)$ and $\Phi_{\boldsymbol{m}}(\boldsymbol{U}^2)$, respectively, then we have
\begin{equation}\label{itfmfgsouniq5}
\begin{split}
\Phi_{\boldsymbol{m}}(\boldsymbol{U}^2)\leq&\sum_{i, j\in\mathcal{S}}c_{ij}(\boldsymbol{m}, \hat{\boldsymbol{P}}^1)m_i\hat{P}^1_{ij}+\sum_{i, j\in\mathcal{S}}U^2_jm_i\hat{P}^1_{ij} \\
                                                                                            =&\sum_{i, j\in\mathcal{S}}c_{ij}(\boldsymbol{m}, \hat{\boldsymbol{P}}^1)m_i\hat{P}^1_{ij}+\sum_{i, j\in\mathcal{S}}U^1_jm_i\hat{P}^1_{ij}+\sum_{i, j\in\mathcal{S}}\left(U^2_j-U^1_j\right)m_i\hat{P}^1_{ij} \\
                                                                                            =&\Phi_{\boldsymbol{m}}(\boldsymbol{U}^1) +\left(\boldsymbol{U}^2-\boldsymbol{U}^1\right)\cdot\Theta_{\boldsymbol{U}^1}(\boldsymbol{m}).
\end{split}
\end{equation}
$\hfill\square$

In addition to Lemma 4.1, we consider the following assumption, which also appears in the reference \cite{GOMES2010308}.

\textbf{(A4)} For any $\boldsymbol{U}^1, \boldsymbol{U}^2\in\mathbb{R}^s$ and any $\boldsymbol{m}^1, \boldsymbol{m}^2\in\mathbb{P}$, there exists a constant $\gamma>0$ that holds the following inequality 
\begin{equation}\label{itfmfgsouniq6}
\boldsymbol{m}^1\cdot\left(\Gamma_{\boldsymbol{m}^1}(\boldsymbol{U}^2)-\Gamma_{\boldsymbol{m}^2}(\boldsymbol{U}^2)\right)+\boldsymbol{m}^2\cdot\left(\Gamma_{\boldsymbol{m}^2}(\boldsymbol{U}^1)-\Gamma_{\boldsymbol{m}^1}(\boldsymbol{U}^1)\right)\geq\gamma\Vert\boldsymbol{m}^1-\boldsymbol{m}^2\Vert^2.
\end{equation}
Then, the uniqueness theorem for the solution to the social optima problem of finite MFGs with given initial-terminal values is as follows.

\textbf{Theorem 4.2.} Suppose (A3) and (A4) hold, given any initial distribution $\boldsymbol{m}^0=\tilde{\boldsymbol{m}}\in\mathbb{P}$ and terminal cost $\boldsymbol{U}^N=\boldsymbol{G}^N\in\mathbb{R}^s$, if $\left(\hat{\boldsymbol{m}}^{n, 1}, \hat{\boldsymbol{U}}^{n, 1}\right)_{n\in\mathcal{N}\cup\{N\}}$ and $\left(\hat{\boldsymbol{m}}^{n, 2}, \hat{\boldsymbol{U}}^{n, 2}\right)_{n\in\mathcal{N}\cup\{N\}}$ are solutions to the social optima problem of finite MFGs with given initial-terminal values (problem (P1)), then for any $n\in\mathcal{N}\cup\{N\}$, $\hat{\boldsymbol{m}}^{n, 1}=\hat{\boldsymbol{m}}^{n, 2}$, $\hat{\boldsymbol{U}}^{n, 1}=\hat{\boldsymbol{U}}^{n, 2}$.

\textbf{Proof.} For symbol simplification, let $\hat{H}^n=\mathop{min}_{\boldsymbol{P}^n}H^n(\boldsymbol{m}^n, \boldsymbol{P}^n, \boldsymbol{U}^{n+1})$. if $\left(\hat{\boldsymbol{m}}^{n, 1}, \hat{\boldsymbol{U}}^{n, 1}\right)_{n\in\mathcal{N}\cup\{N\}}$ and $\left(\hat{\boldsymbol{m}}^{n, 2}, \hat{\boldsymbol{U}}^{n, 2}\right)_{n\in\mathcal{N}\cup\{N\}}$ are solutions to (P1), we have
\begin{equation}\nonumber
\begin{split}
\hat{H}^{n, 1}=\Phi_{\hat{\boldsymbol{m}}^{n, 1}}(\hat{\boldsymbol{U}}^{n+1, 1}), \quad\quad \hat{\boldsymbol{m}}^{n+1, 1}=\Theta_{\hat{\boldsymbol{U}}^{n+1, 1}}(\hat{\boldsymbol{m}}^{n, 1}),  \\
\hat{H}^{n, 2}=\Phi_{\hat{\boldsymbol{m}}^{n, 2}}(\hat{\boldsymbol{U}}^{n+1, 2}), \quad\quad \hat{\boldsymbol{m}}^{n+1, 2}=\Theta_{\hat{\boldsymbol{U}}^{n+1, 2}}(\hat{\boldsymbol{m}}^{n, 2}).
\end{split}
\end{equation}
Then
\begin{equation}\label{itfmfgsouniq7}
\begin{split}
&\!\sum^{N-1}_{n=0}\!\left(\hat{\boldsymbol{U}}^{n+1, 1}\!\!-\!\hat{\boldsymbol{U}}^{n+1, 2}\right)\!\!\left[\left(\Theta_{\hat{\boldsymbol{U}}^{n+1, 1}}(\hat{\boldsymbol{m}}^{n, 1})\!-\!\hat{\boldsymbol{m}}^{n+1, 1}\right)\!\!-\!\!\left(\Theta_{\hat{\boldsymbol{U}}^{n+1, 2}}(\hat{\boldsymbol{m}}^{n, 2})\!-\!\hat{\boldsymbol{m}}^{n+1, 2}\right)\right]\!\!+ \\
&\!\sum^{N-1}_{n=0}\!\left[\left(\Phi_{\hat{\boldsymbol{m}}^{n, 1}}(\hat{\boldsymbol{U}}^{n+1, 1})-\hat{H}^{n, 1}\right)-\left(\Phi_{\hat{\boldsymbol{m}}^{n, 2}}(\hat{\boldsymbol{U}}^{n+1, 2})-\hat{H}^{n, 2}\right)\right]+  \\
&\!\sum^{N-1}_{n=0}\!\left[\left(\Phi_{\hat{\boldsymbol{m}}^{n, 2}}(\hat{\boldsymbol{U}}^{n+1, 2})-\hat{H}^{n, 2}\right)-\left(\Phi_{\hat{\boldsymbol{m}}^{n, 1}}(\hat{\boldsymbol{U}}^{n+1, 1})-\hat{H}^{n, 1}\right)\right]=0.
\end{split}
\end{equation}
Through some mathematical calculations, it can be obtained that
\begin{equation}\label{itfmfgsouniq8}
\begin{split}
&\sum^{N-1}_{n=0}\left[\Phi_{\hat{\boldsymbol{m}}^{n, 1}}(\hat{\boldsymbol{U}}^{n+1, 2})\!-\!\Phi_{\hat{\boldsymbol{m}}^{n, 1}}(\hat{\boldsymbol{U}}^{n+1, 1})\!-\!\left(\hat{\boldsymbol{U}}^{n+1, 2}\!-\!\hat{\boldsymbol{U}}^{n+1, 1}\right)\!\cdot\!\Theta_{\hat{\boldsymbol{U}}^{n+1, 1}}(\hat{\boldsymbol{m}}^{n, 1})\right]\!+ \\
&\sum^{N-1}_{n=0}\left[\Phi_{\hat{\boldsymbol{m}}^{n, 2}}(\hat{\boldsymbol{U}}^{n+1, 1})\!-\!\Phi_{\hat{\boldsymbol{m}}^{n, 2}}(\hat{\boldsymbol{U}}^{n+1, 2})\!-\!\left(\hat{\boldsymbol{U}}^{n+1, 1}\!-\!\hat{\boldsymbol{U}}^{n+1, 2}\right)\!\cdot\!\Theta_{\hat{\boldsymbol{U}}^{n+1, 2}}(\hat{\boldsymbol{m}}^{n, 2})\right]\!+ \\
&\sum^{N-1}_{n=0}\left[\Phi_{\hat{\boldsymbol{m}}^{n, 1}}(\hat{\boldsymbol{U}}^{n+1, 1})\!-\!\Phi_{\hat{\boldsymbol{m}}^{n, 2}}(\hat{\boldsymbol{U}}^{n+1, 1})\!+\!\Phi_{\hat{\boldsymbol{m}}^{n, 2}}(\hat{\boldsymbol{U}}^{n+1, 2})\!-\!\Phi_{\hat{\boldsymbol{m}}^{n, 1}}(\hat{\boldsymbol{U}}^{n+1, 2})-\right. \\
&\left.(\hat{\boldsymbol{m}}^{n+1, 1}-\hat{\boldsymbol{m}}^{n+1, 2})\cdot\left(\hat{\boldsymbol{U}}^{n+1, 1}-\hat{\boldsymbol{U}}^{n+1, 2}\right)\right]=0.
\end{split}
\end{equation}
According to Lemma 4.1,
\begin{equation}\nonumber
\begin{split}
\Phi_{\hat{\boldsymbol{m}}^{n, 1}}(\hat{\boldsymbol{U}}^{n+1, 2})-\Phi_{\hat{\boldsymbol{m}}^{n, 1}}(\hat{\boldsymbol{U}}^{n+1, 1})-\left(\hat{\boldsymbol{U}}^{n+1, 2}-\hat{\boldsymbol{U}}^{n+1, 1}\right)\cdot\Theta_{\hat{\boldsymbol{U}}^{n+1, 1}}(\hat{\boldsymbol{m}}^{n, 1})\leq0, \\
\Phi_{\hat{\boldsymbol{m}}^{n, 2}}(\hat{\boldsymbol{U}}^{n+1, 1})-\Phi_{\hat{\boldsymbol{m}}^{n, 2}}(\hat{\boldsymbol{U}}^{n+1, 2})-\left(\hat{\boldsymbol{U}}^{n+1, 1}-\hat{\boldsymbol{U}}^{n+1, 2}\right)\cdot\Theta_{\hat{\boldsymbol{U}}^{n+1, 2}}(\hat{\boldsymbol{m}}^{n, 2})\leq0.
\end{split}
\end{equation}
Then
\begin{equation}\label{itfmfgsouniq9}
\begin{split}
&\sum^{N-1}_{n=0}\left[\Phi_{\hat{\boldsymbol{m}}^{n, 1}}(\hat{\boldsymbol{U}}^{n+1, 1})-\Phi_{\hat{\boldsymbol{m}}^{n, 2}}(\hat{\boldsymbol{U}}^{n+1, 1})+\Phi_{\hat{\boldsymbol{m}}^{n, 2}}(\hat{\boldsymbol{U}}^{n+1, 2})-\Phi_{\hat{\boldsymbol{m}}^{n, 1}}(\hat{\boldsymbol{U}}^{n+1, 2})-\right. \\
&\left.(\hat{\boldsymbol{m}}^{n+1, 1}-\hat{\boldsymbol{m}}^{n+1, 2})\cdot\left(\hat{\boldsymbol{U}}^{n+1, 1}-\hat{\boldsymbol{U}}^{n+1, 2}\right)\right]\geq0.
\end{split}
\end{equation}
Due to the initial state distribution and the terminal cost are given, i.e. $\hat{\boldsymbol{m}}^{0, 1}=\hat{\boldsymbol{m}}^{0, 2}=\tilde{\boldsymbol{m}}$, $\hat{\boldsymbol{U}}^{N, 1}=\hat{\boldsymbol{U}}^{N, 2}=\boldsymbol{G}^N$, then
\begin{equation}\nonumber
\Phi_{\hat{\boldsymbol{m}}^{0, 1}}(\hat{\boldsymbol{U}}^{1, 1})-\Phi_{\hat{\boldsymbol{m}}^{0, 2}}(\hat{\boldsymbol{U}}^{1, 1})+\Phi_{\hat{\boldsymbol{m}}^{0, 2}}(\hat{\boldsymbol{U}}^{1, 2})-\Phi_{\hat{\boldsymbol{m}}^{0, 1}}(\hat{\boldsymbol{U}}^{1, 2})=0,
\end{equation}
\begin{equation}\nonumber
(\hat{\boldsymbol{m}}^{N, 1}-\hat{\boldsymbol{m}}^{N, 2})\cdot\left(\hat{\boldsymbol{U}}^{N, 1}-\hat{\boldsymbol{U}}^{N, 2}\right)=0.
\end{equation}
From this, it can be concluded that
\begin{equation}\label{itfmfgsouniq10}
\begin{split}
  &\sum^{N-1}_{n=0}\left[\Phi_{\hat{\boldsymbol{m}}^{n, 1}}(\hat{\boldsymbol{U}}^{n+1, 1})-\Phi_{\hat{\boldsymbol{m}}^{n, 2}}(\hat{\boldsymbol{U}}^{n+1, 1})+\Phi_{\hat{\boldsymbol{m}}^{n, 2}}(\hat{\boldsymbol{U}}^{n+1, 2})-\Phi_{\hat{\boldsymbol{m}}^{n, 1}}(\hat{\boldsymbol{U}}^{n+1, 2})-\right. \\
  &\left.(\hat{\boldsymbol{m}}^{n+1, 1}-\hat{\boldsymbol{m}}^{n+1, 2})\cdot\left(\hat{\boldsymbol{U}}^{n+1, 1}-\hat{\boldsymbol{U}}^{n+1, 2}\right)\right] \\
=&\sum^{N-1}_{n=1}\left[\hat{H}^{n, 1}\!-\!\Phi_{\hat{\boldsymbol{m}}^{n, 2}}(\hat{\boldsymbol{U}}^{n+1, 1})\!+\!\hat{H}^{n, 2}\!-\!\Phi_{\hat{\boldsymbol{m}}^{n, 1}}(\hat{\boldsymbol{U}}^{n+1, 2})\!-\!\hat{H}^{n, 1}-\hat{H}^{n, 2}+\right. \\
  &\quad\quad\left.\hat{\boldsymbol{m}}^{n, 1}\cdot\hat{\boldsymbol{U}}^{n, 2}+\hat{\boldsymbol{m}}^{n, 2}\cdot\hat{\boldsymbol{U}}^{n, 1}\right] \\
=&\sum^{N-1}_{n=1}\left[-\hat{\boldsymbol{m}}^{n, 2}\!\cdot\!\Gamma_{\hat{\boldsymbol{m}}^{n, 2}}(\hat{\boldsymbol{U}}^{n+1, 1})\!-\!\hat{\boldsymbol{m}}^{n, 1}\!\cdot\!\Gamma_{\hat{\boldsymbol{m}}^{n, 1}}(\hat{\boldsymbol{U}}^{n+1, 2})\!+\!\hat{\boldsymbol{m}}^{n, 1}\!\cdot\!\Gamma_{\hat{\boldsymbol{m}}^{n, 2}}(\hat{\boldsymbol{U}}^{n+1, 2})+\right. \\
  &\quad\quad\left.\hat{\boldsymbol{m}}^{n, 2}\!\cdot\!\Gamma_{\hat{\boldsymbol{m}}^{n, 1}}(\hat{\boldsymbol{U}}^{n+1, 1})\right]  \\
=&\sum^{N-1}_{n=1}\left[\hat{\boldsymbol{m}}^{n, 1}\cdot\left(\Gamma_{\hat{\boldsymbol{m}}^{n, 2}}(\hat{\boldsymbol{U}}^{n+1, 2})-\Gamma_{\hat{\boldsymbol{m}}^{n, 1}}(\hat{\boldsymbol{U}}^{n+1, 2})\right)+\hat{\boldsymbol{m}}^{n, 2}\cdot\left(\Gamma_{\hat{\boldsymbol{m}}^{n, 1}}(\hat{\boldsymbol{U}}^{n+1, 1})-\right.\right. \\
&\quad\quad\left.\left. \Gamma_{\hat{\boldsymbol{m}}^{n, 2}}(\hat{\boldsymbol{U}}^{n+1, 1})\right)\right] \\
\geq&0.
\end{split}
\end{equation}
Then according to \eqref{itfmfgsouniq6} in assumption (A4)
\begin{equation}\label{itfmfgsouniq11}
\sum^{N-1}_{n=1}-\gamma\Vert\hat{\boldsymbol{m}}^{n, 1}-\hat{\boldsymbol{m}}^{n, 2}\Vert^2\geq0,
\end{equation}
this implies that for any $n\in\mathcal{N}\cup\{N\}$, there is $\hat{\boldsymbol{m}}^{n, 1}=\hat{\boldsymbol{m}}^{n, 2}$. Then iterating backwards from $\hat{\boldsymbol{U}}^{N, 1}=\hat{\boldsymbol{U}}^{N, 2}=\boldsymbol{G}^N$ using $\hat{H}^{n, 1}=\Phi_{\hat{\boldsymbol{m}}^{n, 1}}(\hat{\boldsymbol{U}}^{n+1, 1})$ and $\hat{H}^{n, 2}=\Phi_{\hat{\boldsymbol{m}}^{n, 2}}(\hat{\boldsymbol{U}}^{n+1, 2})$, it can been obtain that for any $n\in\mathcal{N}\cup\{N\}$, the optimal strategies $\hat{\boldsymbol{P}}^{n, 1}=\hat{\boldsymbol{P}}^{n, 2}$, and the individual cost $\hat{\boldsymbol{U}}^{n, 1}=\hat{\boldsymbol{U}}^{n, 2}$.
$\hfill\square$

\section{Stationary solutions to the social optima problem of finite MFGs}
\label{section5}
In this section, we study the existence and uniqueness of stationary solutions to the social optima problem of finite MFGs. Similar to the previous section, compared to the competitive model \cite{GOMES2010308}, the uniqueness condition for stationary solutions corresponding to social optimum be extended from the transition cost $c_{i*}(\boldsymbol{m}^n, \boldsymbol{P}^n_i)$ is only related to $\boldsymbol{P}^n_i$ to being related to $\boldsymbol{P}^n$. For the stationary solution, we adopt a definition similar to \cite{GOMES2010308}. That is, the distribution of each state no longer changes over time, and the individual cost of each state increases by the same constant per time step.

\textbf{Definition 5.1.} A pair of vectors $(\bar{\boldsymbol{m}},\bar{\boldsymbol{U}})$ is a stationary solution to the social optima problem of finite MFGs if there exists a constant $\bar{\lambda}$, called critical value, such that
\begin{equation}\label{fmfgsoss1}
\left\{
\begin{split}
&\bar{\boldsymbol{U}}+\bar{\boldsymbol{\lambda}}=\Gamma_{\bar{\boldsymbol{m}}}(\bar{\boldsymbol{U}}), \\
&\bar{\boldsymbol{m}}=\Theta_{\bar{\boldsymbol{U}}}(\bar{\boldsymbol{m}}),
\end{split}
\right.
\end{equation}
where $\bar{\boldsymbol{\lambda}}=[\bar{\lambda},...,\bar{\lambda}]\in\mathbb{R}^s$, the optimal strategy $\bar{\boldsymbol{P}}=\mathop{argmin}_{\boldsymbol{P}}H(\bar{\boldsymbol{m}}, \boldsymbol{P}, \bar{\boldsymbol{U}})$.

The difference between Definition 5.1 and the definition in reference \cite{GOMES2010308} is that the optimal strategy adopted in Definition 5.1 corresponds to the optimal social cost. For any constant $a\in\mathbb{R}$, the following equation holds
\begin{equation}\label{fmfgsoss2}
\begin{split}
\Phi_{\boldsymbol{m}}(\boldsymbol{U}+\boldsymbol{a})&=\mathop{min}_{\boldsymbol{P}\in\mathbb{P}^s}\sum_{i,j\in\mathcal{S}}c_{ij}(\boldsymbol{m}, \boldsymbol{P})m_iP_{ij}+\sum_{i,j\in\mathcal{S}}(U_j+a)m_iP_{ij} \\
                                                                                 &=\mathop{min}_{\boldsymbol{P}\in\mathbb{P}^s}\left[\sum_{i,j\in\mathcal{S}}c_{ij}(\boldsymbol{m}, \boldsymbol{P})m_iP_{ij}+\sum_{i,j\in\mathcal{S}}U_jm_iP_{ij}\right]+a \\
                                                                                 &=\Phi_{\boldsymbol{m}}(\boldsymbol{U})+a,
\end{split}
\end{equation}
where $\boldsymbol{a}=[a,...,a]\in\mathbb{R}^s$. By utilizing the optimal strategy in \eqref{fmfgsoss2}, it can be inferred that
\begin{equation}\label{fmfgsoss3}
\Gamma_{\boldsymbol{m}}(\boldsymbol{U}+\boldsymbol{a})=\Gamma_{\boldsymbol{m}}(\boldsymbol{U})+\boldsymbol{a}.
\end{equation}
Define an equivalence relation $C$ on $\mathbb{R}^s$: given any constant $a$, $xCy$ represents that the components of $x$ and $y$ differ by the same integer multiple of $a$. Based on this equivalence relation, $\mathbb{R}^s$ can be divided into different equivalence classes. Let $\mathbb{R}^s/\mathbb{R}$ represent the quotient set of these equivalent classes. Then the first equation in \eqref{fmfgsoss1} can be written in $\mathbb{R}^s/\mathbb{R}$ as $\bar{\boldsymbol{U}}=\Gamma_{\bar{\boldsymbol{m}}}(\bar{\boldsymbol{U}})$. This is similar to reference \cite{GOMES2010308}.

\subsection{The critical value}
\label{section5.1}
The following proposition provides a representation of the critical value in Definition 5.1.

\textbf{Proposition 5.1.} If the assumption (A3) holds, and $(\bar{\boldsymbol{m}}, \bar{\boldsymbol{U}})$ is a stationary solution of \eqref{fmfgsoss1}, then the critical value $\bar{\lambda}$ can be expressed as 
\begin{equation}\nonumber
\bar{\lambda}=\sum_{i,j\in\mathcal{S}}c_{ij}(\bar{\boldsymbol{m}},\bar{\boldsymbol{P}})\bar{m}_i\bar{P}_{ij},
\end{equation}
where $\bar{\boldsymbol{P}}$ is the optimal strategy corresponding to social optimum.

\textbf{Proof.} If $(\bar{\boldsymbol{m}}, \bar{\boldsymbol{U}})$ is a stationary solution of \eqref{fmfgsoss1}, then for any $i\in\mathcal{S}$,
\begin{equation}\label{fmfgsoss4}
\bar{U}_i+\bar{\lambda}=\Gamma_{\bar{\boldsymbol{m}}}(\bar{\boldsymbol{U}})_i=\sum_{j\in\mathcal{S}}c_{ij}(\bar{\boldsymbol{m}},\bar{\boldsymbol{P}})\bar{P}_{ij}+\sum_{j\in\mathcal{S}}\bar{U}_j\bar{P}_{ij}.
\end{equation}
Multiplying both sides of the above equation by $m_i$ and adding, we have
\begin{equation}\label{fmfgsoss5}
\sum_{i\in\mathcal{S}}\bar{m}_i\bar{U}_i+\sum_{i\in\mathcal{S}}\bar{m}_i\bar{\lambda}=\sum_{i,j\in\mathcal{S}}c_{ij}(\bar{\boldsymbol{m}},\bar{\boldsymbol{P}})\bar{m}_i\bar{P}_{ij}+\sum_{i, j\in\mathcal{S}}\bar{U}_j\bar{m}_i\bar{P}_{ij}.
\end{equation}
Let $\bar{H}=\sum_{i\in\mathcal{S}}\bar{m}_i\bar{U}_i$, the above equation can be rewritten as
\begin{equation}\label{fmfgsoss6}
\begin{split}
\bar{H}+\bar{\lambda}&=\sum_{i,j\in\mathcal{S}}c_{ij}(\bar{\boldsymbol{m}},\bar{\boldsymbol{U}})\bar{m}_i\bar{P}_{ij}+\sum_{j\in\mathcal{S}}\bar{U}_j\sum_{i\in\mathcal{S}}\bar{m}_i\bar{P}_{ij} \\
                        &=\sum_{i,j\in\mathcal{S}}c_{ij}(\bar{\boldsymbol{m}},\bar{\boldsymbol{U}})\bar{m}_i\bar{P}_{ij}+\bar{H}.
\end{split}
\end{equation}
Therefore,
\begin{equation}\label{fmfgsoss7}
\bar{\lambda}=\sum_{i,j\in\mathcal{S}}c_{ij}(\bar{\boldsymbol{m}},\bar{\boldsymbol{P}})\bar{m}_i\bar{P}_{ij},
\end{equation}
$\hfill\square$

\subsection{Existence of stationary solutions}
\label{section5.2}
The assumptions required in this section are listed below, with the second one being the same as reference \cite{GOMES2010308}.

\textbf{(A5)} For any $\boldsymbol{m}\in\mathbb{P}$, $\boldsymbol{P}\in\mathbb{P}^s$, $i,j\in\mathcal{S}$, $c_{ij}(\boldsymbol{m}, \boldsymbol{P})$ is $C^1$ continuity with respect to $\boldsymbol{P}$.

\textbf{(A6)} For any $\boldsymbol{m}\in\mathbb{P}$, $\boldsymbol{P}\in\mathbb{P}^s$, replace the $p'$-th row of $\boldsymbol{P}$ with its $p$-th row, while keeping the remaining rows (including $p$-th row) unchanged. Set the transformed matrix as $\tilde{\boldsymbol{P}}$. There exists a constant $C>0$, and for any $p,p'\in\mathcal{S}$, the following inequality holds
\begin{equation}\label{fmfgsoss8}
\sum_{j\in\mathcal{S}}\left|c_{pj}(\boldsymbol{m}, \boldsymbol{P})-c_{p'j}(\boldsymbol{m}, \tilde{\boldsymbol{P}})\right|P_{pj}\leq C.
\end{equation}

\textbf{Lemma 5.1.} If $\boldsymbol{m}\in\mathbb{P}$ satisfies assumptions (A3), (A5), (A6), and a constant $\varepsilon$ exists, which satisfies $0<\varepsilon\leq\mathop{min}\{m_i\}_{i\in\mathcal{S}}$, then for any $\boldsymbol{U}\in\mathbb{R}^s$, $p, p'\in\mathcal{S}$, there exists a constant $C'>0$, such that the following inequality holds
\begin{equation}\label{fmfgsoss9}
\left|\Gamma_{\boldsymbol{m}}(\boldsymbol{U})_p-\Gamma_{\boldsymbol{m}}(\boldsymbol{U})_{p'}\right|\leq C'.
\end{equation}

\textbf{Proof.} If $\boldsymbol{m}\in\mathbb{P}$ satisfies assumptions (A3),  there exists a unique optimal strategy $\boldsymbol{P}$ corresponding to the social optimum. For any $p, p'\in\mathcal{S}$, let $\tilde{\boldsymbol{P}}$ represent a state transition matrix, whose $p'$-th row is the $p$-th row of $\boldsymbol{P}$, and the other rows (including $p$-th row) are the same as $\boldsymbol{P}$, then we have
\begin{equation}\label{fmfgsoss10}
\sum_{i, j\in\mathcal{S}}\!c_{ij}(\boldsymbol{m}, \boldsymbol{P})m_iP_{ij}\!+\!\sum_{i, j\in\mathcal{S}}\!U_jm_iP_{ij}\leq\sum_{i, j\in\mathcal{S}}\!c_{ij}(\boldsymbol{m}, \tilde{\boldsymbol{P}})m_i\tilde{P}_{ij}\!+\!\sum_{i, j\in\mathcal{S}}\!U_jm_i\tilde{P}_{ij}.
\end{equation}
Since $\tilde{\boldsymbol{P}}$ and $\boldsymbol{P}$ are identical except for the $p'$-th row, it can be inferred from the above equation
\begin{equation}\label{fmfgsoss11}
\begin{split}
&\sum_{(i\neq p'), j\in\mathcal{S}}\!\!\!\!\left(c_{ij}(\boldsymbol{m}, \boldsymbol{P})\!-\!c_{ij}(\boldsymbol{m}, \tilde{\boldsymbol{P}})\right)\!m_iP_{ij}\!+\!\!\sum_{j\in\mathcal{S}}\!\left(c_{p'j}(\boldsymbol{m}, \boldsymbol{P})P_{p'j}\!-\!c_{p'j}(\boldsymbol{m}, \tilde{\boldsymbol{P}})\tilde{P}_{p'j}\right)\!m_{p'} \\
&+\sum_{j\in\mathcal{S}}U_jm_{p'}P_{p'j}-\sum_{j\in\mathcal{S}}U_jm_{p'}\tilde{P}_{p'j}\leq0.
\end{split}
\end{equation}
According to assumptions (A5), for any $i, j\in\mathcal{S}$, $c_{ij}(\boldsymbol{m}, \boldsymbol{P})$ is first order continuously differentiable with respect to $\boldsymbol{P}$. And because the feasible set of the strategy is a bounded closed set, there exists a Lipschitz constant $L>0$, such that
\begin{equation}\label{fmfgsoss12}
\vert c_{ij}(\boldsymbol{m}, \boldsymbol{P})-c_{ij}(\boldsymbol{m}, \tilde{\boldsymbol{P}})\vert\leq L\Vert\boldsymbol{P}-\tilde{\boldsymbol{P}}\Vert=L\sqrt{\sum_{j\in\mathcal{S}}(P_{p'j}-\tilde{P}_{p'j})^2}\leq L\sqrt{s}.
\end{equation}
Take the largest one of the Lipschitz constants corresponding to any $i\neq p', j\in\mathcal{S}$, and still record it as $L$, then
\begin{equation}\label{fmfgsoss13}
\begin{split}
&-L\sqrt{s}\left(\sum_{(i\neq p'), j\in\mathcal{S}}m_iP_{ij}\right)+\sum_{j\in\mathcal{S}}c_{p'j}(\boldsymbol{m}, \boldsymbol{P})m_{p'}P_{p'j}+\sum_{j\in\mathcal{S}}U_jm_{p'}P_{p'j}- \\
&\left(\sum_{j\in\mathcal{S}}c_{p'j}(\boldsymbol{m}, \tilde{\boldsymbol{P}})m_{p'}\tilde{P}_{p'j}+\sum_{j\in\mathcal{S}}U_jm_{p'}\tilde{P}_{p'j}\right)\leq0.
\end{split}
\end{equation}
Since there exists a constant $\varepsilon$, which satisfies $0<\varepsilon\leq\mathop{min}\{m_i\}_{i\in\mathcal{S}}$, and the $p'$-th row of $\tilde{\boldsymbol{P}}$ is the same as the $p$-th row of $\boldsymbol{P}$, then we can deduce that
\begin{equation}\label{fmfgsoss14}
\Gamma_{\boldsymbol{m}}(\boldsymbol{U})_{p'}\leq\sum_{j\in\mathcal{S}}c_{p'j}(\boldsymbol{m}, \tilde{\boldsymbol{P}})P_{pj}+\sum_{j\in\mathcal{S}}U_jP_{pj}+\frac{L\sqrt{s}}{\varepsilon}.
\end{equation}
Hence
\begin{equation}\label{fmfgsoss15}
\Gamma_{\boldsymbol{m}}(\boldsymbol{U})_{p'}-\Gamma_{\boldsymbol{m}}(\boldsymbol{U})_p\leq\sum_{j\in\mathcal{S}}\left[c_{p'j}(\boldsymbol{m}, \tilde{\boldsymbol{P}})-c_{pj}(\boldsymbol{m}, \boldsymbol{P})\right]P_{pj}+\frac{L\sqrt{s}}{\varepsilon}\leq C+\frac{L\sqrt{s}}{\varepsilon}.
\end{equation}
Based on the assumption (A6), the second less-than or equal sign in the above equation holds. Let $C'=C+(L\sqrt{s})/\varepsilon$, exchange the role of $p$ and $p'$, and after synthesis, we can obtain, for any $p, p'\in\mathcal{S}$,
\begin{equation}\label{fmfgsoss16}
\left|\Gamma_{\boldsymbol{m}}(\boldsymbol{U})_p-\Gamma_{\boldsymbol{m}}(\boldsymbol{U})_{p'}\right|\leq C'.
\end{equation}
$\hfill\square$

According to Lemma 5.1, we can use a method similar to reference \cite{GOMES2010308} to prove the existence of stationary solutions based on Brower's fixed point theorem.  The following theorem provides a sufficient condition for the existence of stationary solutions to the social optima problem of finite MFGs. 

\textbf{Theorem 5.1.} If the assumptions (A2), (A3), (A5), and (A6) hold, there exists a pair of vectors $(\bar{\boldsymbol{m}}, \bar{\boldsymbol{U}})$ that holds \eqref{fmfgsoss1}, i.e. it is a stationary solution to the social optima problem of finite MFGs. 

\textbf{Proof.} The proof process can refer to Theorem 3 in reference \cite{GOMES2010308}.
$\hfill\square$

\subsection{Uniqueness of stationary solutions}
\label{section5.3}
Before proving the uniqueness of stationary solutions to the social optima problem of finite MFGs,  Lemma 4.1 needs to be investigated again. For \eqref{itfmfgsouniq4}, we cannot strictly write it as
\begin{equation}\nonumber
\Phi_{\boldsymbol{m}}(\boldsymbol{U}^2)-\Phi_{\boldsymbol{m}}(\boldsymbol{U}^1)-\left(\boldsymbol{U}^2-\boldsymbol{U}^1\right)\cdot\Theta_{\boldsymbol{U}^1}(\boldsymbol{m})\leq-\eta\Vert\boldsymbol{U}^2-\boldsymbol{U}^1\Vert^2,
\end{equation}
where $\eta>0$. This is mainly due to the characteristics of $\Phi_{\boldsymbol{m}}(\cdot)$ with respect to $V+constant$ shown in \eqref{fmfgsoss2}. In fact, set $\boldsymbol{U}^2=\boldsymbol{U}^1+\boldsymbol{k}$, where $k$ is a constant and $\boldsymbol{k}=[k, ..., k]\in\mathbb{R}^s$, then using \eqref{fmfgsoss2} we can deduce that
\begin{equation}\nonumber
\Phi_{\boldsymbol{m}}(\boldsymbol{U}^2)-\Phi_{\boldsymbol{m}}(\boldsymbol{U}^1)-\left(\boldsymbol{U}^2-\boldsymbol{U}^1\right)\cdot\Theta_{\boldsymbol{U}^1}(\boldsymbol{m})=0.
\end{equation}
Define the following norm on the quotient set $\mathbb{R}^s/\mathbb{R}$ composed of equivalent classes mentioned at the beginning of this section
\begin{equation}\label{fmfgsoss17}
\Vert\boldsymbol{U}\Vert_\#:=\mathop{inf}_{\boldsymbol{\nu}\in\mathbb{R}^s}\Vert\boldsymbol{U}+\boldsymbol{\nu}\Vert,
\end{equation}
where $\boldsymbol{\nu}=[\nu, ... ,\nu]\in\mathbb{R}^s$. However, similar to reference \cite{GOMES2010308}, we can make the following stricter assumption

\textbf{(A7)} For any $\boldsymbol{m}\in\mathbb{P}$, $\boldsymbol{U}^1, \boldsymbol{U}^2\in\mathbb{R}^s$, the following inequality holds
\begin{equation}\nonumber
\Phi_{\boldsymbol{m}}(\boldsymbol{U}^2)-\Phi_{\boldsymbol{m}}(\boldsymbol{U}^1)-\left(\boldsymbol{U}^2-\boldsymbol{U}^1\right)\cdot\Theta_{\boldsymbol{U}^1}(\boldsymbol{m})\leq-\eta\Vert\boldsymbol{U}^2-\boldsymbol{U}^1\Vert^2_\#,
\end{equation}
where $\eta>0$.

\textbf{Theorem 5.2.} Suppose assumption (A3), (A4), and (A7) hold, if both vector pairs $(\bar{\boldsymbol{m}}^1$, $\bar{\boldsymbol{U}}^1), (\bar{\boldsymbol{m}}^2, \bar{\boldsymbol{U}}^2)$ are stationary solutions, i.e. they satisfy
\begin{equation}\nonumber
\begin{split}
\bar{\boldsymbol{U}}^1+\bar{\boldsymbol{\lambda}}^1=\Gamma_{\bar{\boldsymbol{m}}^1}(\bar{\boldsymbol{U}}^1), \quad\quad \bar{\boldsymbol{m}}^1=\Theta_{\bar{\boldsymbol{U}}^1}(\bar{\boldsymbol{m}}^1),  \\
\bar{\boldsymbol{U}}^2+\bar{\boldsymbol{\lambda}}^2=\Gamma_{\bar{\boldsymbol{m}}^2}(\bar{\boldsymbol{U}}^2), \quad\quad \bar{\boldsymbol{m}}^2=\Theta_{\bar{\boldsymbol{U}}^2}(\bar{\boldsymbol{m}}^2),
\end{split}
\end{equation}
where the critical values $\bar{\lambda}^1, \bar{\lambda}^2$ are constants. Then $\bar{\boldsymbol{m}}^1=\bar{\boldsymbol{m}}^2$, $\bar{\lambda}^1=\bar{\lambda}^2$, $\bar{\boldsymbol{U}}^2=\bar{\boldsymbol{U}}^1+\boldsymbol{k}$, $k$ is a contant, $\boldsymbol{k}=[k, ..., k]\in\mathbb{R}^s$.

\textbf{Proof.} \textbf{(1) Prove that $\bar{\boldsymbol{m}}^1=\bar{\boldsymbol{m}}^2$.} Let $\bar{H}^1=\sum_{i\in\mathcal{S}}\bar{m}^1_i\bar{U}^1_i$, $\bar{H}^2=\sum_{i\in\mathcal{S}}\bar{m}^2_i\bar{U}^2_i$, according to \eqref{fmfgsoss5}, we have
\begin{equation}\nonumber
\bar{H}^1+\bar{\lambda}^1=\Phi_{\bar{\boldsymbol{m}}^1}(\bar{\boldsymbol{U}}^1), \quad\quad \bar{H}^2+\bar{\lambda}^2=\Phi_{\bar{\boldsymbol{m}}^2}(\bar{\boldsymbol{U}}^2).
\end{equation}
Hence
\begin{equation}\label{fmfgsoss18}
\begin{split}
&(\bar{\boldsymbol{U}}^1-\bar{\boldsymbol{U}}^2)\cdot\left[\left(\Theta_{\bar{\boldsymbol{U}}^1}(\bar{\boldsymbol{m}}^1)- \bar{\boldsymbol{m}}^1\right)-\left(\Theta_{\bar{\boldsymbol{U}}^2}(\bar{\boldsymbol{m}}^2)-\bar{\boldsymbol{m}}^2\right)\right]+\left[\left(\Phi_{\bar{\boldsymbol{m}}^1}(\bar{\boldsymbol{U}}^1)-\bar{H}^1\right)-\right. \\
&\left.\left(\Phi_{\bar{\boldsymbol{m}}^2}(\bar{\boldsymbol{U}}^2)-\bar{H}^2\right)\right]+\left[\left(\Phi_{\bar{\boldsymbol{m}}^2}(\bar{\boldsymbol{U}}^2)-\bar{H}^2\right)-\left(\Phi_{\bar{\boldsymbol{m}}^1}(\bar{\boldsymbol{U}}^1)-\bar{H}^1\right)\right]=0.
\end{split}
\end{equation}
Through some mathematical calculations, it can be obtained that
\begin{equation}\label{fmfgsoss19}
\begin{split}
&\left[\Phi_{\bar{\boldsymbol{m}}^1}(\bar{\boldsymbol{U}}^2)-\Phi_{\bar{\boldsymbol{m}}^1}(\bar{\boldsymbol{U}}^1)-(\bar{\boldsymbol{U}}^2-\bar{\boldsymbol{U}}^1)\cdot\Theta_{\bar{\boldsymbol{U}}^1}(\bar{\boldsymbol{m}}^1)\right]+ \\
&\left[\Phi_{\bar{\boldsymbol{m}}^2}(\bar{\boldsymbol{U}}^1)-\Phi_{\bar{\boldsymbol{m}}^2}(\bar{\boldsymbol{U}}^2)-(\bar{\boldsymbol{U}}^1-\bar{\boldsymbol{U}}^2)\cdot\Theta_{\bar{\boldsymbol{U}}^2}(\bar{\boldsymbol{m}}^2)\right]+ \\
&\left[\Phi_{\bar{\boldsymbol{m}}^1}(\bar{\boldsymbol{U}}^1)-\Phi_{\bar{\boldsymbol{m}}^2}(\bar{\boldsymbol{U}}^1)+\Phi_{\bar{\boldsymbol{m}}^2}(\bar{\boldsymbol{U}}^2)-\Phi_{\bar{\boldsymbol{m}}^1}(\bar{\boldsymbol{U}}^2)\!-\!(\bar{\boldsymbol{m}}^1\!-\!\bar{\boldsymbol{m}}^2)\cdot(\bar{\boldsymbol{U}}^1\!-\!\bar{\boldsymbol{U}}^2)\right]\!=\!0.
\end{split}
\end{equation}
According to Lemma 4.1, the first rows in \eqref{fmfgsoss19} are both less than or equal to zero. Then
\begin{equation}\label{fmfgsoss20}
\Phi_{\bar{\boldsymbol{m}}^1}(\bar{\boldsymbol{U}}^1)-\Phi_{\bar{\boldsymbol{m}}^2}(\bar{\boldsymbol{U}}^1)+\Phi_{\bar{\boldsymbol{m}}^2}(\bar{\boldsymbol{U}}^2)-\Phi_{\bar{\boldsymbol{m}}^1}(\bar{\boldsymbol{U}}^2)-(\bar{\boldsymbol{m}}^1-\bar{\boldsymbol{m}}^2)\cdot(\bar{\boldsymbol{U}}^1-\bar{\boldsymbol{U}}^2)\geq0.
\end{equation}
Since
\begin{equation}\label{fmfgsoss21}
\begin{split}
&\Phi_{\bar{\boldsymbol{m}}^1}(\bar{\boldsymbol{U}}^1)-\Phi_{\bar{\boldsymbol{m}}^2}(\bar{\boldsymbol{U}}^1)+\Phi_{\bar{\boldsymbol{m}}^2}(\bar{\boldsymbol{U}}^2)-\Phi_{\bar{\boldsymbol{m}}^1}(\bar{\boldsymbol{U}}^2)-(\bar{\boldsymbol{m}}^1-\bar{\boldsymbol{m}}^2)\cdot(\bar{\boldsymbol{U}}^1-\bar{\boldsymbol{U}}^2) \\
=&\bar{H}^1+\bar{\lambda}^1+\bar{H}^2+\bar{\lambda}^2-\Phi_{\bar{\boldsymbol{m}}^1}(\bar{\boldsymbol{U}}^2)-\Phi_{\bar{\boldsymbol{m}}^2}(\bar{\boldsymbol{U}}^1)-\bar{H}^1-\bar{H}^2+\\
   &\bar{\boldsymbol{m}}^1\cdot(\Gamma_{\bar{\boldsymbol{m}}^2}(\bar{\boldsymbol{U}}^2)-\bar{\boldsymbol{\lambda}}^2)+\bar{\boldsymbol{m}}^2\cdot(\Gamma_{\bar{\boldsymbol{m}}^1}(\bar{\boldsymbol{U}}^1)-\bar{\boldsymbol{\lambda}}^1)  \\ 
=&\bar{\boldsymbol{m}}^1\cdot\Gamma_{\bar{\boldsymbol{m}}^2}(\bar{\boldsymbol{U}}^2)+\bar{\boldsymbol{m}}^2\cdot\Gamma_{\bar{\boldsymbol{m}}^1}(\bar{\boldsymbol{U}}^1)-\Phi_{\bar{\boldsymbol{m}}^1}(\bar{\boldsymbol{U}}^2)-\Phi_{\bar{\boldsymbol{m}}^2}(\bar{\boldsymbol{U}}^1) \\
=&\bar{\boldsymbol{m}}^1\cdot\left(\Gamma_{\bar{\boldsymbol{m}}^2}(\bar{\boldsymbol{U}}^2)-\Gamma_{\bar{\boldsymbol{m}}^1}(\bar{\boldsymbol{U}}^2)\right)+\bar{\boldsymbol{m}}^2\cdot\left(\Gamma_{\bar{\boldsymbol{m}}^1}(\bar{\boldsymbol{U}}^1)-\Gamma_{\bar{\boldsymbol{m}}^2}(\bar{\boldsymbol{U}}^1)\right) \\
\geq&0,
\end{split}
\end{equation}
according to assumption (A4), it can be concluded that
\begin{equation}\label{fmfgsoss22}
-\gamma\Vert\bar{\boldsymbol{m}}^1-\bar{\boldsymbol{m}}^2\Vert^2\geq0.
\end{equation}
Since $\gamma$ in the above equation is a constant greater than zero, then $\bar{\boldsymbol{m}}^1=\bar{\boldsymbol{m}}^2$.

\textbf{(2) Prove that $\bar{\lambda}^1=\bar{\lambda}^2$.} Set $\bar{\boldsymbol{m}}=\bar{\boldsymbol{m}}^1=\bar{\boldsymbol{m}}^2$, $\bar{\boldsymbol{P}}^2$ is the optimal strategy corresponding to $\Phi_{\bar{\boldsymbol{m}}}(\bar{\boldsymbol{U}}^2)$, then according to Lemma 4.1, we have
\begin{equation}\label{fmfgsoss23}
\begin{split}
  &\Phi_{\bar{\boldsymbol{m}}}(\bar{\boldsymbol{U}}^1)-\Phi_{\bar{\boldsymbol{m}}}(\bar{\boldsymbol{U}}^2)-\sum_{i, j\in\mathcal{S}}(\bar{U}^1_j-\bar{U}^2_j)\bar{m}_i\bar{P}^2_{ij} \\
=&\left(\Phi_{\bar{\boldsymbol{m}}}(\bar{\boldsymbol{U}}^1)-\sum_{j\in\mathcal{S}}\bar{m}_j\bar{U}^1_j\right)-\left(\Phi_{\bar{\boldsymbol{m}}}(\bar{\boldsymbol{U}}^2)-\sum_{j\in\mathcal{S}}\bar{m}_j\bar{U}^2_j\right) \\
=&\bar{\lambda}^1-\bar{\lambda}^2 \\
\leq&0.
\end{split}
\end{equation}
Similarly, set $\bar{\boldsymbol{P}}^1$ represents the optimal strategy corresponding to $\Phi_{\bar{\boldsymbol{m}}}(\bar{\boldsymbol{U}}^1)$, we have
\begin{equation}\label{fmfgsoss24}
\begin{split}
  &\Phi_{\bar{\boldsymbol{m}}}(\bar{\boldsymbol{U}}^2)-\Phi_{\bar{\boldsymbol{m}}}(\bar{\boldsymbol{U}}^1)-\sum_{i, j\in\mathcal{S}}(\bar{U}^2_j-\bar{U}^1_j)\bar{m}_i\bar{P}^1_{ij} \\
=&\left(\Phi_{\bar{\boldsymbol{m}}}(\bar{\boldsymbol{U}}^2)-\sum_{j\in\mathcal{S}}\bar{m}_j\bar{U}^2_j\right)-\left(\Phi_{\bar{\boldsymbol{m}}}(\bar{\boldsymbol{U}}^1)-\sum_{j\in\mathcal{S}}\bar{m}_j\bar{U}^1_j\right) \\
=&\bar{\lambda}^2-\bar{\lambda}^1 \\
\leq&0.
\end{split}
\end{equation}
Then $\bar{\lambda}^1=\bar{\lambda}^2$.

\textbf{(3) Prove that $\bar{\boldsymbol{U}}^2=\bar{\boldsymbol{U}}^1+\boldsymbol{k}$.} Set $\bar{\boldsymbol{m}}=\bar{\boldsymbol{m}}^1=\bar{\boldsymbol{m}}^2$, based on \eqref{itfmfgsouniq6}, \eqref{fmfgsoss19}, and \eqref{fmfgsoss21}, we can obtain
\begin{equation}\label{fmfgsoss25}
\begin{split}
&\left[\Phi_{\bar{\boldsymbol{m}}}(\bar{\boldsymbol{U}}^2)-\Phi_{\bar{\boldsymbol{m}}}(\bar{\boldsymbol{U}}^1)-(\bar{\boldsymbol{U}}^2-\bar{\boldsymbol{U}}^1)\cdot\Theta_{\bar{\boldsymbol{U}}^1}(\bar{\boldsymbol{m}})\right]+\left[\Phi_{\bar{\boldsymbol{m}}}(\bar{\boldsymbol{U}}^1)-\Phi_{\bar{\boldsymbol{m}}}(\bar{\boldsymbol{U}}^2)-\right. \\
&\left.(\bar{\boldsymbol{U}}^1-\bar{\boldsymbol{U}}^2)\cdot\Theta_{\bar{\boldsymbol{U}}^2}(\bar{\boldsymbol{m}})\right]=0.
\end{split}
\end{equation}
Then according to the assumption (A7),
\begin{equation}\label{fmfgsoss26}
-\eta^1\Vert\bar{\boldsymbol{U}}^2-\bar{\boldsymbol{U}}^1\Vert^2_\#-\eta^2\Vert\bar{\boldsymbol{U}}^1-\bar{\boldsymbol{U}}^2\Vert^2_\#\geq0.
\end{equation}
Since $\eta^1, \eta^2>0$, then $\Vert\bar{\boldsymbol{U}}^2-\bar{\boldsymbol{U}}^1\Vert^2_\#=\Vert\bar{\boldsymbol{U}}^1-\bar{\boldsymbol{U}}^2\Vert^2_\#=0$. In view of the definition of the norm $\Vert\cdot\Vert_\#$, it can be obtained that there exists a constant $k$, such that $\bar{\boldsymbol{U}}^2=\bar{\boldsymbol{U}}^1+\boldsymbol{k}$.
$\hfill\square$

\section{Examples}
\label{section6}
This section gives two examples of social optima problems of finite MFGs. Example 1 corresponds to the existence and uniqueness of the solutions satisfying the given initial-terminal value. Example 2 corresponds to the existence and uniqueness of stationary solutions.
\subsection{Example 1}
\label{section6.1}
Set $\mathcal{S}=\{1, 2, 3\}$. For any $i, j\in\mathcal{S}$, let
\begin{equation}\label{exone1}
c_{ij}(\boldsymbol{m}, \boldsymbol{P})=\alpha_1m_i+\alpha_2(\sum_{p\in\mathcal{S}}m_pP_{pj})-\alpha_3\frac{P_{ij}}{lnP_{ij}},
\end{equation}
where $\alpha_i>0, i\in\mathcal{S}$ are the weight values. It is easy to know that $c_{ij}(\boldsymbol{m}, \boldsymbol{P})$ is a differentiable convex function with respect to $\boldsymbol{P}$. For any $\boldsymbol{P}^1, \boldsymbol{P}^2\in\mathbb{P}^3, \boldsymbol{P}^1\neq\boldsymbol{P}^2$,
\begin{equation} \label{exone2}
\begin{split}
   &\sum_{i, j\in\mathcal{S}}\left(c_{ij}(\boldsymbol{m}, \boldsymbol{P}^1)-c_{ij}(\boldsymbol{m}, \boldsymbol{P}^2)\right)\left(P^1_{ij}-P^2_{ij}\right)m_i \\
=&\sum_{i, j\in\mathcal{S}}\alpha_2\left(\sum_{p\in\mathcal{S}}m_pP_{pj}^1-\sum_{p\in\mathcal{S}}m_pP_{pj}^2\right)\left(P^1_{ij}-P^2_{ij}\right)m_i+ \\
   &\sum_{i, j\in\mathcal{S}}\alpha_3\left(-\frac{P^1_{ij}}{lnP^1_{ij}}+\frac{P^2_{ij}}{lnP^2_{ij}}\right)\left(P^1_{ij}-P^2_{ij}\right)m_i \\
=&\sum_{j\in\mathcal{S}}\alpha_2\left(\sum_{p\in\mathcal{S}}m_pP_{pj}^1\!-\!\sum_{p\in\mathcal{S}}m_pP_{pj}^2\right)^2\!\!+\!\!\sum_{i, j\in\mathcal{S}}\alpha_3\left(-\frac{P^1_{ij}}{lnP^1_{ij}}\!+\!\frac{P^2_{ij}}{lnP^2_{ij}}\right)\left(P^1_{ij}\!-\!P^2_{ij}\right)m_i \\
>&0.
\end{split}
\end{equation}
Since $\boldsymbol{P}^1\neq\boldsymbol{P}^2$, and the function $f(x)=-x/lnx$ is a strictly monotonically increasing function when $x\in(0, 1)$, the greater than sign in the above equation holds. According to the function setting of $c_{ij}(\boldsymbol{m}, \boldsymbol{P})$, $P_{ij}>0$. Then take a sufficiently small $\varepsilon$, such that $P_{ij}\geq\varepsilon$. Furthermore, considering the constraint $\sum_{j\in\mathcal{S}}P_{ij}=1$, it is known that there is an optimal strategy $\boldsymbol{P}$ in the feasible set. Similar to the uniqueness proof in Theorem 3.1, it can be seen that the optimal strategy $\boldsymbol{P}$ is unique.

If the initial state distribution and the terminal cost are given, for any $n\in\mathcal{N}$, we have 
\begin{equation}\label{exone3}
\begin{split}
H^n(\boldsymbol{m}^n, \boldsymbol{P}^n, \boldsymbol{U}^{n+1})=&\sum_{i, j\in\mathcal{S}}\left(\alpha_1m^n_i+\alpha_2\left(\sum_{p\in\mathcal{S}}m^n_pP^n_{pj}\right)-\alpha_3\frac{P^n_{ij}}{lnP^n_{ij}}\right)m^n_iP^n_{ij}+\\
&\sum_{i, j\in\mathcal{S}}U^{n+1}_jm^n_iP^n_{ij}.
\end{split}
\end{equation}
This is a continuous function on $\mathbb{P}\times\mathbb{P}^s\times\mathbb{R}^s$. According to Theorem 4.1, for the social optima problem of the finite mean field with given initial-terminal values, there exists a sequence of pairs of vectors $\left(\hat{\boldsymbol{m}}^n, \hat{\boldsymbol{U}}^n\right)_{n\in\mathcal{N}\cup\{N\}}$ that is its solution.

If \eqref{exone1} is changed as follows
\begin{equation}\label{exone4}
c_{ij}(\boldsymbol{m}, \boldsymbol{P})=\alpha_1m_i-\alpha_2\frac{P_{ij}}{lnP_{ij}},
\end{equation}
we can prove the uniqueness of the solution to the given initial-terminal values problem mentioned above. In this case 
\begin{equation}\label{exone5}
\Phi_{\boldsymbol{m}}(\boldsymbol{U})=\sum_{i\in\mathcal{S}}\alpha_1(m_i)^2+\mathop{min}_{\boldsymbol{P}}\sum_{i, j\in\mathcal{S}}\left(-\alpha_2\frac{P_{ij}}{lnP_{ij}}+U_j\right)m_iP_{ij}.
\end{equation}
By using the KKT condition to solve, it can be seen that the value of the optimal strategy $\hat{\boldsymbol{P}}$ is only related to $\boldsymbol{U}$. So for $\boldsymbol{m}^1$ and $\boldsymbol{m}^2$ that belong to the feasible set, and for any $\boldsymbol{U}^1, \boldsymbol{U}^2\in\mathbb{R}^3$, if the constant $\gamma\leq\alpha_1$, we have
\begin{equation}\label{exone6}
\begin{split}
   &\boldsymbol{m}^1\cdot\left(\Gamma_{\boldsymbol{m}^1}(\boldsymbol{U}^2)-\Gamma_{\boldsymbol{m}^2}(\boldsymbol{U}^2)\right)+\boldsymbol{m}^2\cdot\left(\Gamma_{\boldsymbol{m}^2}(\boldsymbol{U}^1)-\Gamma_{\boldsymbol{m}^1}(\boldsymbol{U}^1)\right)  \\
=&\sum_{i\in\mathcal{S}}\alpha_1m^1_i(m^1_i-m^2_i)+\sum_{i\in\mathcal{S}}\alpha_1m^2_i(m^2_i-m^1_i) \\
\geq&\gamma\Vert\boldsymbol{m}^1-\boldsymbol{m}^2\Vert^2.
\end{split}
\end{equation}
Then according to Theorem 4.2, the solution of the social optima problem of the finite mean field with given initial-terminal values is unique.

\subsection{Example 2}
\label{section6.2}
Set $\mathcal{S}=\{1, 2\}$. For any $i, j\in\mathcal{S}$, let
\begin{equation}\label{extwo1}
c_{ij}(\boldsymbol{m}, \boldsymbol{P})=m_i+\left(\sum_{(p\neq i)\in\mathcal{S}}m_pP_{pj}\right)+ln(P_{ij})^2.
\end{equation}
Then at time $n$,
\begin{equation}\label{extwo2}
\begin{split}
H^n(\boldsymbol{m}^n, \boldsymbol{P}^n, \boldsymbol{U}^{n+1})=&\sum_{i, j\in\mathcal{S}}\left(m^n_i+\left(\sum_{(p\neq i)\in\mathcal{S}}m^n_pP^n_{pj}\right)+ln(P^n_{ij})^2\right)m^n_iP^n_{ij}+\\
&\sum_{i, j\in\mathcal{S}}U^{n+1}_jm^n_iP^n_{ij}.
\end{split}
\end{equation}
This is a continuous function on $\mathbb{P}\times\mathbb{P}^s\times\mathbb{R}^s$. Through some mathematical calculations, it is easy to know that $H^n(\boldsymbol{m}^n, \boldsymbol{P}^n, \boldsymbol{U}^{n+1})$ is a strictly convex function with respect to $\boldsymbol{P}^n$. We construct a bounded closed convex feasible set similar to Example 1, then  there exists a unique optimal strategy $\hat{\boldsymbol{P}}^n$ at time $n$. In fact, the optimal strategy is given by
\begin{equation}\label{extwo3}
\hat{\boldsymbol{P}}^n_{ij}=\frac{e^{-\left(\sum_{(p\neq i)\in\mathcal{S}}m^n_pP^n_{pj}\right)-\frac{U^{n+1}_j}{2}}}{\sum_{q\in\mathcal{S}}e^{-\left(\sum_{(p\neq i)\in\mathcal{S}}m^n_pP^n_{pq}\right)-\frac{U^{n+1}_q}{2}}}.
\end{equation}

Next, we prove the existence of stationary solutions. Obviously, $c_{ij}(\boldsymbol{m}, \boldsymbol{P})$ is $C^1$ continuity with respect to $\boldsymbol{P}$. Let $\hat{\boldsymbol{P}}$ be the optimal strategy corresponding to the social optima at the current time. For any $p, p'\in\mathcal{S}$, let $\tilde{\boldsymbol{P}}$ represent a state transition matrix, whose $p'$-th row is the $p$-th row of $\hat{\boldsymbol{P}}$, and the other rows (including $p$-th row) are the same as $\hat{\boldsymbol{P}}$, then we have
\begin{equation}\label{extwo4}
\begin{split}
   &\sum_{j\in\mathcal{S}}\left|c_{pj}(\boldsymbol{m}, \hat{\boldsymbol{P}})-c_{p'j}(\boldsymbol{m}, \tilde{\boldsymbol{P}})\right|\hat{P}_{pj} \\
=&\sum_{j\in\mathcal{S}}\left|m_p+\left(\sum_{(k\neq p)\in\mathcal{S}}m_k\hat{P}_{kj}\right)-m_{p'}-\left(\sum_{(k\neq p')\in\mathcal{S}}m_k\tilde{P}_{kj}\right)\right|\hat{P}_{pj} \\
=&\sum_{j\in\mathcal{S}}\left|m_p(1-\hat{P}_{pj})-m_{p'}(1-\hat{P}_{p'j})\right|\hat{P}_{pj} \\
\leq&C.
\end{split}
\end{equation}
Since $m_p, m_{p'}, \hat{P}_{pj}\in(0, 1)$, it can be inferred that the constant $C$ in the above inequality exists. According to Lemma 5.1 and Theorem 5.1, there exists a pair of vectors $(\bar{\boldsymbol{m}}, \bar{\boldsymbol{U}})$ is a stationary solution.

If \eqref{extwo1} is changed as follows
\begin{equation}\label{extwo5}
c_{ij}(\boldsymbol{m}, \boldsymbol{P})=m_i+ln(P_{ij})^2,
\end{equation}
then similar to Example 1, we can obtain that the assumption (A3) holds, similar to Section 4.4 of the reference \cite{GOMES2010308}, we can obtain that the assumption (A6) holds. According to Theorem 5.2, the stationary solution is unique.

\section{Conclusions}

This paper discusses the social optimum of finite mean field games. A sufficient condition for the existence and uniqueness of the optimal strategy for minimizing the social cost is provided. The existence and uniqueness conditions of equilibrium solutions under initial-terminal constraints in the finite horizon and stationary solutions in the infinite horizon are given. Finally, two examples corresponding respectively to equilibrium solutions and stationary solutions are provided.



\bibliographystyle{elsarticle-num} 
\bibliography{references}





\end{document}